\documentclass[12pt]{amsart}
\usepackage{amssymb,latexsym,cite}
\usepackage[usenames, dvipsnames]{color}
\definecolor{mygray}{gray}{0.6}

\usepackage{amssymb,amsfonts,amsthm, amsmath, breqn, color, float, mathtools, booktabs, hvfloat, url,cite}

\advance\textheight by \topskip   \numberwithin{equation}{section}

\newtheorem{theorem}{Theorem}[section]

\newtheorem{lemma}[theorem]{Lemma}

\newcommand{\beq}{\begin{eqnarray*}}
\newcommand{\feq}{\end{eqnarray*}}
\newcommand{\beqn}{\begin{eqnarray}}
\newcommand{\feqn}{\end{eqnarray}}

\def\obrace{\iftrue{\else}\fi}
\def\cbrace{\iffalse{\else}\fi}

\let\originalparagraph\paragraph
\renewcommand{\paragraph}[2][.]{\originalparagraph{#2#1}}

\newcommand{\zz}{{\mathbb Z}}

\newcommand{\cc}{{\mathbb C}}

\newcommand{\bes}{\begin{split}}
\newcommand{\fes}{\end{split}}

\newtheorem*{conj*}{Conjecture}
\makeatletter \@addtoreset{theorem}{section}\makeatother

\makeatletter \@addtoreset{theorem}{section}\makeatother

\makeatletter \@addtoreset{theorem}{section}\makeatother

\newtheorem*{theorem*}{Theorem}

\makeatletter
\newcommand{\leqnomode}{\tagsleft@true\let\veqno\@@leqno}
\newcommand{\reqnomode}{\tagsleft@false\let\veqno\@@eqno}
\makeatother
\usepackage{algorithm}
\usepackage{graphicx}
\usepackage[noend]{algpseudocode}
\makeatletter
\def\BState{\State\hskip-\ALG@thistlm}
\makeatother

\DeclareCaptionLabelFormat{noname}{#2}
\newlength\myindent
\def\CP{\mathcal{CP}}
\begin{document}
\title[Convex polyominoes]{Convex polyominoes revisited: Enumeration of outer site perimeter, interior vertices, and boundary vertices of certain degrees}

\author{Toufik Mansour}
\address{Department of Mathematics, University of Haifa, 3498838 Haifa, Israel}
\email{tmansour@univ.haifa.ac.il }

\author{Reza Rastegar}
\address{Occidental Petroleum Corporation, Houston, TX 77046 and Departments of Mathematics and Engineering, University of Tulsa, OK 74104, USA - Adjunct Professor}
\email{reza\_rastegar2@oxy.com}


\begin{abstract}
The main contribution of this paper is a new column-by-column method for the decomposition of generating functions of convex polyominoes suitable for enumeration with respect to various statistics including but not limited to interior vertices, boundary vertices of certain degrees, and outer site perimeter. Using this decomposition, among other things, we show that

A) the average number of interior vertices
over all convex polyominoes of perimeter $2n$ is asymptotic
to $\frac{n^2}{12}+\frac{n\sqrt{n}}{3\sqrt{\pi}}
-\frac{(21\pi-16)n}{12\pi}.$

B) the average number of boundary vertices with degree two over all convex polyominoes of perimeter $2n$ is asymptotic to
$\frac{n+6}{2}+\frac{1}{\sqrt{\pi n}}+\frac{(16-7\pi)}{4\pi n}.$
Additionally, we obtain an explicit generating function counting the number of convex polyominoes with $n$ boundary vertices of degrees at most three and show that this number is asymptotic to
$
\frac{n+1}{40}\left(\frac{3+\sqrt{5}}{2}\right)^{n-3}
+\frac{\sqrt[4]{5}(2-\sqrt{5})}{80\sqrt{\pi n}}\left(\frac{3+\sqrt{5}}{2}\right)^{n-2}.
$
Moreover, we show that the expected number of the boundary vertices of degree four over all convex polyominoes with $n$  vertices of degrees at most three is asymptotically
$
\frac{n}{\sqrt{5}}-\frac{\sqrt[4]{125}(\sqrt{5}-1)\sqrt{n}}{10\sqrt{\pi}}.
$

C) the number of convex polyominoes with the outer-site perimeter $n$ is asymptotic to
$\frac{3(\sqrt{5}-1)}{20\sqrt{\pi
		n}\sqrt[4]{5}}\left(\frac{3+\sqrt{5}}{2}\right)^n,$
and show the expected number of the outer-site perimeter over all convex polyominoes with perimeter $2n$ is asymptotic to $\frac{25n}{16}+\frac{\sqrt{n}}{4\sqrt{\pi}}+\frac{1}{8}.$
Lastly, we prove that the expected perimeter over all convex polyominoes with the outer-site perimeter $n$ is asymptotic
to $\sqrt[4]{5}n$.	
\end{abstract}
\maketitle

\section{Introduction}

Interest in general squared lattice polyominoes, also known as polyominoes or animals, arose around the same time in statistical physics and combinatorics communities. In the former, they are mainly used for modeling percolation  \cite{bro1}, the mechanics of macro-molecules \cite{Temperly} such as the collapse of branched polymers \cite{peard1}, and the analysis of growth models \cite{eden1, klar1, Read}. In combinatorics, they initially appeared as puzzles and mathematical games \cite{GO}, however, soon they became about the enumeration of graphs \cite{har1}. For the earliest works and some history, see \cite{melou4, gut9, Viennot, conway1} and references therein.  Among many mathematical mysteries about polyominoes, a long-standing and seems-to-be-intractable open question is to determine the number of general polyominoes of exactly $n$ cells, whose growth rate, though proven to exist in \cite{klar1}, is still unknown. To tackle this long-lasting challenge, various constraints such as directional convexity and/or directional growth have been deployed to reduce the difficulty of the enumeration process in the general polyominoes. In this paper we study an important subset of polyominoes with such constraints, namely {\it convex polyominoes,}  that are among  one of the most interesting and studied subclasses of polyominoes satisfying several convexity constraints. In particular, we enumerate them with respect to {\it interior vertices}, {\it vertices of various degrees} on the boundary, and {\it outer-site perimeter} and the interplay with {\it perimeter}. \\ \par

In what follows we first review a non-exhaustive selection of some important results on this subject. In one of the first pioneering work,  Bender \cite{ben1} showed that the number of convex polyominoes with exactly $n$ cells is of asymptotic form $\alpha \beta^n$, where $\alpha \approx 2.67564$ and $\beta \approx 2.30914.$ Later, Delest and Viennot \cite{DV} used a bijection between convex polyominoes and words of an algebraic language and showed the generating function for convex polyominoes indexed by the half-perimeter is
\beqn \label{cp_perim_gf}
\frac{x^2(1 - 8x + 21x^2 - 19x^3 + 4x^4)}{(1 - 2x)(1 - 4x)^2} - \frac{2x^4}{(1-4x)\sqrt{1-4x}}.
\feqn
This expression was obtained by differencing two series with positive
terms, whose combinatorial interpretation was given by Bousquet-M\'elou and
Guttmann in \cite{melou3}. The closed formula for the number of convex polyominoes with perimeter $2n+8$ is therefore given by
\beqn \label{eqfcpf}
(2n+11)4^n-4(2n+1){2n \choose n}.
\feqn
In addition, Bousquet-M\'elou derived several expressions for the generating functions of convex polyominoes according to the area, the number of rows and columns \cite{melou1, melou2}.  We remark that \eqref{eqfcpf} has been also derived using other different methods. For example, Kim \cite{kim1} provided an elementary proof of this formula by counting the pairs of non-crossing lattice paths in plane. In \cite{ler1}, the authors used Burnside's lemma and counted the number of unique convex polyominoes of $n$ cells up to both  reflections and rotations. In \cite{lin1}, Lin and Chang provided a generating function for the number of convex polyominoes bounded by a given rectangle. Recently, Buchin et al. \cite{buchin1} extended their work and gave a new combinatorial proof. They also counted the subclass of these polyominoes that contain the lower left corner of the enclosing rectangle (directed polyominoes) and calculated the first and second moments of the number of common points of two monotone lattice paths between two given points. Del Lungo et al. \cite{lungo1} defined the generating function of convex polyominoes according to the semi-perimeter using the ECO method \cite{eco, EB}. Finally, Hochst\"attler et al. \cite{lo1} provided an efficient method to generate convex polyominoes at random in polynomial time.\\
\par

The main contribution of this paper is a new column-by-column method for decomposition of generating functions of convex polyominoes with respect to various statistics listed in Table \ref{tab1}. Layer-by-layer based decomposition methods have been used previously in the context of enumeration of various combinatorial objects including but not limited to polyominoes \cite{doron1, doron2, melou1}. The success and effectiveness of each of these decompositions are highly dependent on the enumeration problem and its specific structure. For instance, using our method we are able to track several two-vertex statistics simultaneously and in a more compact manner, which seems to provide stronger results comparing with other known decompositions. In order to elaborate our methodology, we include four different applications in which we enumerate convex polyominoes with respect to perimeter, interior vertices, and lastly but most importantly, the number of vertices of certain degrees on the boundary and outer-site perimeter. We choose to include the already-well-known first application to elaborate on the methodology through a simple example and also collect several useful intermediate results for the subsequent applications. A few definitions and various discussions about the results are in order. \\
\par
\begin{figure}[t]
	\begin{picture}(20,25)
	\put(0,5){
		\multiput(0,0)(0,4){2}{\put(0,0){\line(1,0){4}}\put(4,0){\line(0,1){4}}
			\put(4,4){\line(-1,0){4}}\put(0,4){\line(0,-1){4}}}
		\multiput(4,-8)(0,4){5}{\put(0,0){\line(1,0){4}}\put(4,0){\line(0,1){4}}
			\put(4,4){\line(-1,0){4}}\put(0,4){\line(0,-1){4}}}
		\multiput(8,-4)(0,4){6}{\put(0,0){\line(1,0){4}}\put(4,0){\line(0,1){4}}
			\put(4,4){\line(-1,0){4}}\put(0,4){\line(0,-1){4}}}
		\multiput(12,-4)(0,4){6}{\put(0,0){\line(1,0){4}}\put(4,0){\line(0,1){4}}
			\put(4,4){\line(-1,0){4}}\put(0,4){\line(0,-1){4}}}
		\multiput(16,4)(0,4){3}{\put(0,0){\line(1,0){4}}\put(4,0){\line(0,1){4}}
			\put(4,4){\line(-1,0){4}}\put(0,4){\line(0,-1){4}}} }
	\end{picture}
	\caption{An example of a convex polyomino with $22$ cells, a horizontal perimeter of $10$, a vertical perimeter of $14$, and an outer-site perimeter of $18$. It also has $11$ interior vertices, $10$ boundary vertices of degree two, and $6$ boundary vertices of degree four.}
	
	\label{fig1}
\end{figure}
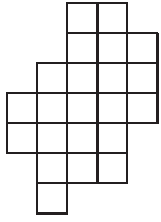

We denote by a {\it cell} $[u, v]$ a unit square in the Cartesian plane $\zz^2$ with its sides parallel to the coordinate axes and with its center at an integer point $(u, v)\in \zz^2$. We refer to the vertices of the square $(u \pm 1/2,v \pm 1/2)$ and its vertical (resp. horizontal) sides as the vertices and vertical (resp. horizontal) edges of cell $[u,v]$. Two cells $[u, v]$ and $[r, s]$ are edge-connected if $|u - r| + |v - s| = 1$. A polyomino is a finite edge-connected set $\nu$ of cells, where for each distinct pair of cells $[u_1,v_1]$ and $[u_2,v_2]$ in $\nu$, there is a finite consecutive sequence of edge-adjacent cells in $\nu$ connecting $[u_1,v_1]$ and $[u_2,v_2]$. A cell in $\nu$ with at least one edge in common with a cell in $\nu^c:=\zz^2 \setminus \nu$ is referred to as a boundary cell of $\nu$. This common edge is referred to as a boundary edge. The degree of a vertex in $\nu$ is the number of edges in $\nu$ incident to that vertex. A vertex in $\nu$ can be of degree $2$, $3$, or $4$. A vertex in a polyomino $\nu$ is called an {\it interior} vertex if it is shared by exactly four different cells of $\nu$, otherwise it is called a {\it boundary} vertex. The {\it perimeter} (resp. {\it semi-perimeter}) of a polyomino $\nu$ is (resp. half of)  the total number of its boundary edges. Similarly, we define the outer-site perimeter of $\nu$ as the number of cells in $\nu^c$ with at least one common edge with a cell in $\nu$. A {\it  column (resp. row)} of a polyomino $\nu$ is a non-empty subset of cells in $\nu$ of from $[c,v]$ (resp. $[u,c]$) for a fixed $c \in \zz$. A polyomino is called {\it column-convex} ({\it row-convex}) if each of its columns (rows) is a single contiguous block of cells. A {\it convex} polyomino is both column-convex and row-convex. It is known that the semi-perimeter of a convex polyomino is the total number of its columns and rows. We use $\CP$ to denote the set of all unique convex polyominoes up-to translation on the squared lattice. See Figure \ref{fig1} for an example of convex polyomino and its corresponding statistics studied (see Table \ref{tab1}) in this paper.

\begin{table}[htp]
\begin{tabular}{c||l}
  Statistic& Description\\\hline\hline
  $a$& number of cells\\
  $c$& number of cells minus $1$ in first column\\
  $h$& number of rows\\
  $v$& number of columns\\
  $int$&number of interior vertices\\
  $d_m$&number vertices of degree $m$ on the boundary\\
  $o$& outer-site perimeter\\\hline
\end{tabular}
\caption{Statistics}\label{tab1}
\end{table}

The rest of the paper is organized as follows. In section \ref{fund_idea}, we provide a few notations and explain the basic idea behind the decomposition.  For our warm-up application, in Section \ref{CP_decom}, we study the generating function
\beq
C(x,y,t;z):=\sum_{\nu\in \CP}x^{v(\nu)}y^{h(\nu)}t^{a(\nu)}z^{c(\nu)},
\feq
where $a(\nu),$ $c(\nu),$ $h(\nu),$ and $v(\nu)$ denote the number of cells, the number of cells minus one in the first column, the number of rows, and the number of columns in the polyomino $\nu$, respectively. This analysis will provide another proof of \eqref{cp_perim_gf} among other things. To study $C(x,y,t;z)$, we introduce two subsets of convex polyominoes; roughly speaking, the sets of convex polyominoes that are staircase-like on either one or both sides. By studying these subsets, we obtain an explicit equation \eqref{eqacp} for $C(x,y,t;z)$ where we study further through the Kernel method \cite{Ban}. 

In Section~\ref{internal_ver}, with a similar decomposition method, we study the generating function
\beq
F(x,y,q;z):=\sum_{\nu\in \CP}x^{v(\nu)}y^{h(\nu)}q^{int(\nu)}z^{c(\nu)},
\feq
where $int(\nu)$ denotes the number of interior vertices of the polyomino $\nu$. We find an exact expression for the generating function $F(x,y,q;z)$ (See Theorem~\ref{internal_thm}), and show that the total number of interior vertices over all convex polyominoes of perimeter $2n$ is
\beqn \label{total_int_exp}
\frac{1}{6}(4n^3-78n^2+77n+321)4^{n-5}+\frac{2}{3}(5n-8)(n-3)\binom{2n-6}{n-3},
\feqn
for all $n\geq5$.  Moreover, by  \eqref{eqfcpf}, we show that the expected number of the interior vertices over all convex polyominoes with perimeter $2n$ is asymptotically
\beq
&& \frac{\frac{1}{6}(4n^3-78n^2+77n+321)4^{n-5}+\frac{2}{3}(5n-8)(n-3)\binom{2n-6}{n-3}}{(2n+3)4^{n-4}-4(2n-7)\binom{2n-8}{n-4}}\notag \\
&& \qquad \sim\frac{n^2}{12}+\frac{n\sqrt{n}}{3\sqrt{\pi}}
-\frac{(21\pi-16)n}{12\pi},
\feq
where $f_n\sim g_n$ if and only if $\lim_{n\rightarrow\infty}f_n/g_n=1$. \\
\par

Recall that understanding of polyominoes with respect to boundary vertices of certain degrees will provide a natural way to study the complexity and shape of a {\it typical} polyomino. For instance, if we traverse the boundary of a polyomino with a large number of boundary vertices of degrees three,  we tend to make less changes in the direction of movement; hence this number can be considered as an index of smoothness of the boundary. Therefore, in Section~\ref{boundary_ver}, we study the generating function that tracks the degree of vertices of the boundary vertices; that is
\beq
D(x,y,q_2,q_3,q_4; z) := \sum_{\nu \in \CP} x^{v(\nu)} y^{h(\nu)} q_2^{d_2(\nu)}
q_3^{d_3(\nu)}q_4^{d_4(\nu)} z^{c(\nu)},
\feq
where $d_m(\nu)$ denotes the number of vertices of degree $m$ on the boundary of $\nu$. We obtain an exact expression for this generating function, and show that the total number of all boundary vertices of degree two of all convex polyominoes with perimeter $2n$ is
\beqn \label{degree2_total}
\frac{1}{2}4^{n-4}(2n^2+15n+11)-\frac{2}{2n-7}(n^2+2n-19)(n-3)\binom{2n-6}{n-3},
\feqn
for all $n\geq5$. Similarly, we prove that the expected number of the boundary vertices of degree two over all convex polyominoes with perimeter $2n$ is given by
\beqn \label{degree2_average}
&&\frac{\frac{1}{2}4^{n-4}(2n^2+15n+11)-\frac{2}{2n-7}(n^2+2n-19)(n-3)\binom{2n-6}{n-3}}{(2n+3)4^{n-4}-4(2n-7)\binom{2n-8}{n-4}}\notag \\
&&\qquad\sim\frac{n+6}{2}+\frac{1}{\sqrt{\pi n}}+\frac{(16-7\pi)}{4\pi n}.
\feqn
Furthermore, we obtain an explicit generating function \eqref{deg3_gf} that counts the number of convex polyominoes with $n$ boundary vertices of degrees at most three. This generating function implies in particular that the number of polyominoes in $\CP$ with $n$ boundary vertices of degrees at most three is asymptotic to
\beqn\label{degree3_asymp}
\frac{n+1}{40}\left(\frac{3+\sqrt{5}}{2}\right)^{n-3}
+\frac{\sqrt[4]{5}(2-\sqrt{5})}{80\sqrt{\pi n}}\left(\frac{3+\sqrt{5}}{2}\right)^{n-2}.
\feqn
Moreover, by considering the generating function $D(1,1,q,q,p;1),$ we show that the number of polyominoes $\nu$ in $\CP$ with $n$ boundary vertices of degrees at most three is equal to the number of polyominoes $\nu$ in $\CP$ with $n-4$ boundary vertices of degrees of minimum three, for all $n\geq4$. In addition, this generating function implies that the total number of all boundary vertices of degrees four over all polyominoes $\nu\in \CP$ with $n$ boundary vertices of degrees at most three is asymptotically
\beqn \label{bv1-asym}
\frac{\sqrt{5}(n+1)(n+2)}{200}\left(\frac{3+\sqrt{5}}{2}\right)^{n-3}
-\frac{\sqrt[4]{125}\sqrt{n^3}}{200\sqrt{\pi}}
\left(\frac{3+\sqrt{5}}{2}\right)^{n-7/2}.
\feqn
This, furthermore, implies that the expected number of the boundary vertices of degree four over all polyominoes in $\CP$ with $n$ boundary vertices of degrees of at most three is asymptotically
\beqn \label{bv2-asym}
\frac{n}{\sqrt{5}}-\frac{\sqrt[4]{125}(\sqrt{5}-1)\sqrt{n}}{10\sqrt{\pi}}.
\feqn \\
\par

For our last application, we enumerate convex polyominoes with respect to their outer-site perimeter and the interplay with perimeter. It appears that the enumeration with respect to outer-site perimeter in convex polyominoe requires tracking the relationship among three columns at the time and hence, in principle, its intrinsic complexity is more challenging comparing with the enumeration with respect to only perimeter. This is also evident by comparing the complexity of enumeration with respect to perimeter for bargraphs. We refer to \cite{BR,DD,Fer2} and references within for a few examples on the enumeration with respect to perimeters in different subsets of polyominoes. To accomplish this enumeration task, in Section \ref{outer_section}, we study the generating function
\beq
J(x,y,q;z):=\sum_{\nu\in \CP} x^{v(\nu)} y^{h(\nu)} q^{o(\nu)} z^{c(\nu)},
\feq
where $o(\nu)$ denotes the outer-site perimeter of the polyomino $\nu$. Using this generating function, in particular, we show that the total sum of the outer-site perimeter over all convex polyominoes with perimeter $2n$ is
\beqn \label{tot_out_resp_perimeter}
(50n^2+79n+105)4^{n-6}+2^{n-6}-\frac{1}{2n-7}(6n^2-19n-8)(n-3).
\feqn
Dividing this by \eqref{eqfcpf}, we then obtain an asymptotic form for the average of the outer-site perimeter over all convex polyominoes of perimeter $2n$; that is,
\beqn \label{avg_out_resp_perimeter}
\frac{25n}{16}+\frac{\sqrt{n}}{4\sqrt{\pi}}+\frac{1}{8}.
\feqn
Similarly, we show that the expected perimeter over all convex polyominoes with outer perimeter $n$ is asymptotic to
\beqn \label{avg_perim_res_outer}
\sqrt[4]{5}n.
\feqn
Lastly, but most importantly, we are also able to find the asymptotic of the number of convex polyominoes with the outer-site perimeter $n$ as
\beqn \label{outer_asympt}
\frac{3(\sqrt{5}-1)}{20\sqrt{\pi
		n}\sqrt[4]{5}}\left(\frac{3+\sqrt{5}}{2}\right)^n.
\feqn

\section{The fundamental idea} \label{fund_idea}

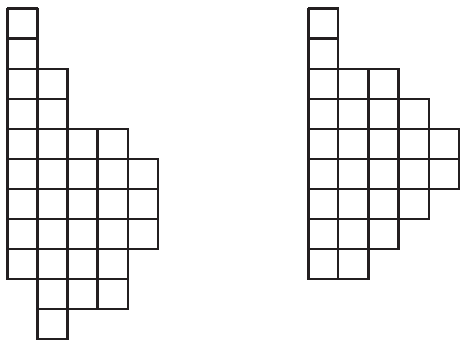
\begin{figure}[htp]
	\begin{picture}(55,50)
	\put(0,5){		\multiput(0,0)(0,4){9}{\put(0,0){\line(1,0){4}}\put(4,0){\line(0,1){4}}
			\put(4,4){\line(-1,0){4}}\put(0,4){\line(0,-1){4}}}		\multiput(4,-8)(0,4){9}{\put(0,0){\line(1,0){4}}\put(4,0){\line(0,1){4}}
			\put(4,4){\line(-1,0){4}}\put(0,4){\line(0,-1){4}}} \multiput(8,-4)(0,4){6}{\put(0,0){\line(1,0){4}}\put(4,0){\line(0,1){4}}
			\put(4,4){\line(-1,0){4}}\put(0,4){\line(0,-1){4}}} \multiput(12,-4)(0,4){6}{\put(0,0){\line(1,0){4}}\put(4,0){\line(0,1){4}}
			\put(4,4){\line(-1,0){4}}\put(0,4){\line(0,-1){4}}} \multiput(16,4)(0,4){3}{\put(0,0){\line(1,0){4}}\put(4,0){\line(0,1){4}}
			\put(4,4){\line(-1,0){4}}\put(0,4){\line(0,-1){4}}} }
\put(40,5){
\multiput(0,0)(0,4){9}{\put(0,0){\line(1,0){4}}\put(4,0){\line(0,1){4}}
			\put(4,4){\line(-1,0){4}}\put(0,4){\line(0,-1){4}}}		\multiput(4,0)(0,4){7}{\put(0,0){\line(1,0){4}}\put(4,0){\line(0,1){4}}
			\put(4,4){\line(-1,0){4}}\put(0,4){\line(0,-1){4}}}
\multiput(8,4)(0,4){6}{\put(0,0){\line(1,0){4}}\put(4,0){\line(0,1){4}}
			\put(4,4){\line(-1,0){4}}\put(0,4){\line(0,-1){4}}} \multiput(12,8)(0,4){4}{\put(0,0){\line(1,0){4}}\put(4,0){\line(0,1){4}}
			\put(4,4){\line(-1,0){4}}\put(0,4){\line(0,-1){4}}} \multiput(16,12)(0,4){2}{\put(0,0){\line(1,0){4}}\put(4,0){\line(0,1){4}}
			\put(4,4){\line(-1,0){4}}\put(0,4){\line(0,-1){4}}} }
	\end{picture}
	\caption{An example of polyominos in (left) $\CP^u$ (right) $\CP^{bu}$  \label{fig10}}
\end{figure}

By the assumption that the polyominoes are unique up-to translation maps, acting on the set $\zz^2,$ we may identify the bottom cell in the first column of all polyominos with the cell $[0,0]$ and say the bottom (resp. top) cell of the $i$-th column of $\nu$ is at the {\it position} $k$ if its bottom (resp. top) edge coincides with the line $y=k+\frac{1}{2}$ (resp. $y=k-\frac{1}{2}$). We denote the position of the bottom (resp. top) cell of the $i$-th column by $b(\nu,i)$ (resp. $u(\nu,i)$). We define $\CP^u$ to be the set of all nonempty convex polyominoes $\nu,$ where for each column $j$ of $\nu$, $u(\nu,s)\leq u(\nu, j)$ for all $s\geq j+1$. Similarly, $\CP^b$ is the set of all nonempty convex polyominoes $\nu$ such that for all columns $j,$ and for all $s\geq j+1,$ $b(\nu,s)\geq b(\nu,j)$. Set $\CP^{bu}:=\CP^u\cap\CP^b$. We note that all $\CP^{u},$ $\CP^{b},$ $\CP^{ub}$ are proper subsets of more known set of stack polyominoes in literatures. See for example \cite{melou1} for a discussion on stack polyominoes.
By a clear upside down symmetry, there is a bijection between the set $\CP^u$ and the set $\CP^b$. See Figure \ref{fig10}, for an example of polyominoes in $\CP^u$ and $\CP^{bu}$.
In addition, we let $\CP_k^u,$ $\CP_k^b,$ $\CP_k^{bu},$ and $\CP_k$ denote the set of polyominoes in $\CP^u,$ $\CP^b,$ $\CP^{bu},$ and $\CP$ with $k$ cells in their first columns. Our enumeration methodology for a given two-vertex statistics works as follows. The first step is to count the number of polyominoes in $\CP^{ub}$ with respect to the statistics of interest. Typically the counting in $\CP^{ub}$ for many two-vertex statistics is straightforward. The second step is to extend the result from $\CP^{ub}$ to $\CP^{u}$ and $\CP^{b}.$ This can be done since the generating functions of the statistics of interest can be written recursively in terms of the corresponding generating functions in $\CP^{ub}$. When this is done, the last step is to obtain the results for $\CP$ by lifting up the result obtained for $\CP^u$ and $\CP^b$. We have chosen four examples to show how these steps are done. The reported results in the last three sections never appeared in literatures before. \par

In each of the following sections, by considering the cases in $\CP^{bu}$, $\CP^{u}$, $\CP^{b}$, $\CP$, we establish a functional system of linear combinations of multivariate series of forms $A(t_1,\ldots, t_r;x_1,\ldots,x_s)$ which do not depends on $x_1,\ldots,x_s,$ simultaneously. The coefficients of each series $A(t_1,\ldots, t_r;x_1,\ldots,x_s)$ are referred to as kernels \cite{doron0}. A systematic approach to solve these class of equations is
given in \cite{Ban}. We remark that throughout this paper, we conduct calculation and manipulations of all generating functions having the following facts in mind. Let $\mathbb{Q}$, $\mathbb{Q}[x_1, \dots, x_s],$ and $\mathbb{Q}[[x_1, \ldots, x_s]]$ denote, respectively, the field of rational numbers, the ring of polynomials in $x_1,\ldots,x_s$, and the ring of formal power
series in $x_1,\ldots,x_s$ with coefficients in $\mathbb{Q}$.
Recall that a series $A\in \mathbb{Q}[x_1, \dots, x_s][[t_1,\ldots, t_r]]$ is $D$-finite if
its partial derivatives span a finite dimensional vector space
over the field of rational functions in $t_1,\cdots t_r$ with coefficients in
$\mathbb{Q}[x_1,\ldots,x_s]$; see \cite{lip1,lip2,St} for an elaborative discussion. Note that any algebraic series is $D$-finite. The specializations of a $D$-finite series, obtained by assigning values in $\mathbb{Q}$ to a subset of variables,
are $D$-finite, if well-defined.  Moreover, if $A$ is $D$-finite, then
any substitution $x_i=1$ with $i\in I\subseteq [s]$ or/and $x_i=x_j$
for $(i,j)\in I\times I\subseteq [s]\times[s]$ into $A$ is also
$D$-finite \cite{lip1}. One last remark is that in some cases, we also perform singularity analysis (see \cite[Section VI]{FS} for a comprehensive review), to derive asymptotic forms for our results. The calculations in these cases are understood in the appropriate domain of complex number $\cc$. We omit the details for the sake of brevity as it is standard. 

\section{Perimeter}
\label{CP_decom}

Our first application is the enumeration of convex polyominoes according to perimeter. The main objective is to elaborate on the methodology and provide several intermediate results that will be also used in three subsequent sections. Additionally, as a side result, we give another proof for \eqref{eqfcpf}.

\subsection{Enumeration in $\CP^{bu}$}
Let $C^{bu}_k(x,y,t) := \sum_{\nu\in \CP_k^{bu}} x^{v(\nu)} y^{h(\nu)} t^{a(\nu)}.$
Define
$$C^{bu}(x,y,t;z):=\sum_{k\geq1}C^{bu}_k(x,y,t)z^{k-1} = \sum_{\nu\in \CP^{bu}} x^{v(\nu)} y^{h(\nu)} t^{a(\nu)} z^{c(\nu)},$$ where $z$ marks the cells in the first column. Consider a polyomino $\nu$ in $\CP^{bu}_k$. It has either one column whose contribution to $C^{bu}_k(x,y,t)$ is exactly $xy^kt^k$, or has more than one column. In the latter case, let $1\leq l \leq k$ be the number of cells in the second column. Since $u(\nu,2)\leq u(\nu,1)$ and $b(\nu,2)\geq b(\nu,1)$, for a given $l$, there are exactly $(k+1-l)$ ways of gluing the second column to the first column. In each of these cases, the contribution of $\nu$ to the generating function $C^{bu}_k(x,y,t)$ is exactly $xy^{k-l}t^kC^{bu}_l(x,y,t).$ Hence,
\beq
C^{bu}_k(x,y,t)=xy^kt^k+\sum_{l=1}^k(k+1-l)xy^{k-l}t^kC^{bu}_l(x,y,t).
\feq
By multiplying by $z^{k-1}$ and summing over $k\geq1$, we obtain
\beqn
C^{bu}(x,y,t;z)&=&\frac{xyt}{1-ytz}+\sum_{k\geq1}\sum_{l=1}^k(k+1-l)xy^{k-l}t^kz^{k-1}C^{bu}_l(x,y,t) \notag\\
&=&\frac{xyt}{1-ytz}+xt\sum_{l\geq1}\sum_{k=l}^\infty (k+1-l)y^{k-l}(tz)^{k-1}C^{bu}_l(x,y,t)\notag \\
&=&\frac{xyt}{1-ytz}+
\frac{xt}{(1-ytz)^2}\sum_{l\geq1}(tz)^{l-1}C^{bu}_l(x,y,t)\notag \\
&=& \frac{xyt}{1-ytz}+ \frac{xt}{(1-ytz)^2}C^{bu}(x,y,t;tz) \label{Cbu_exp}.
\feqn
We therefore obtain the following result
\begin{lemma}\label{lem1}
For $t=1,$
$$C^{bu}(x,y,1;z)=\frac{xy(1-yz)}{(1-yz)^2-x}$$
and for $|t|<1$, we have
$$C^{bu}(x,y,t;z)=\sum_{j\geq0}\frac{y(xt)^{j+1}}{(1-yt^{j+1}z)\prod_{i=1}^{j}(1-yt^iz)^2}.$$
\end{lemma}
Note that the first part is simply obtained by substituting $t=1$ in \eqref{Cbu_exp}. The second part is the result of an iterative application of \eqref{Cbu_exp}  when $|t|<1$. Note  that $C^{bu}(x,x,1;1)=\frac{x^2(1-x)}{1-3x+x^2}$. Hence, as a corollary to this lemma one can re-derive the known-fact that the number of polyominoes in $\CP^{bu}$ with perimeter $2n$ is the $(2n-3)$-st Fibonacci number, where the $n$-th Fibonacci number is defined via the recurrence relation $a_n=a_{n-1}+a_{n-2}$ and with the initial conditions $a_0=0$ and $a_1=1$.

\subsection{Enumeration in $\CP^{u}$}
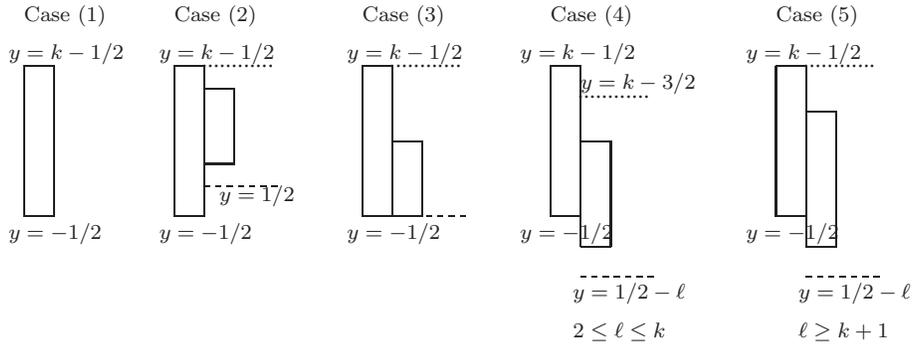
\begin{figure}[htp]
	\begin{picture}(120,45)
	\setlength{\unitlength}{0.1cm}\put(0,1){
		%---1
		\put(0,15){\put(0,26){\tiny  Case (1)}	\put(0,0){\line(1,0){4}}\put(4,0){\line(0,1){20}}\put(4,20){\line(-1,0){4}}\put(0,20){\line(0,-1){20}}
\put(-2,-3){\tiny$y=-1/2$}\put(-2,21){\tiny$y=k-1/2$}}
		%---2
		\put(20,15){\put(0,26){\tiny  Case (2)} \put(0,0){\line(1,0){4}}\put(4,0){\line(0,1){20}}\put(4,20){\line(-1,0){4}}\put(0,20){\line(0,-1){20}}			\put(4,7){\put(0,0){\line(1,0){4}}\put(4,0){\line(0,1){10}}\put(4,10){\line(-1,0){4}}\put(0,10){\line(0,-1){10}}}			\multiput(4,20)(0.8,0){12}{\circle*{.01}}
\multiput(4,4)(1.5,0){7}{\line(1,0){.7}}
\put(6,2){\tiny$y=1/2$}\put(-2,-3){\tiny$y=-1/2$}\put(-2,21){\tiny$y=k-1/2$}}
		%---3
		\put(45,15){\put(0,26){\tiny  Case (3)} \put(0,0){\line(1,0){4}}\put(4,0){\line(0,1){20}}\put(4,20){\line(-1,0){4}}\put(0,20){\line(0,-1){20}}			\put(4,0){\put(0,0){\line(1,0){4}}\put(4,0){\line(0,1){10}}\put(4,10){\line(-1,0){4}}\put(0,10){\line(0,-1){10}}} \multiput(4,20)(0.8,0){12}{\circle*{.01}}
\multiput(4,0)(1.5,0){7}{\line(1,0){.7}}
\put(-2,-3){\tiny$y=-1/2$}\put(-2,21){\tiny $y=k-1/2$}}
		%---4
		\put(70,15){\put(0,26){\tiny  Case (4)} \put(0,0){\line(1,0){4}}\put(4,0){\line(0,1){20}}\put(4,20){\line(-1,0){4}}\put(0,20){\line(0,-1){20}} \put(4,-4){\put(0,0){\line(1,0){4}}\put(4,0){\line(0,1){14}}\put(4,14){\line(-1,0){4}}\put(0,14){\line(0,-1){14}}} \multiput(4,16)(0.8,0){12}{\circle*{.01}}
\multiput(4,-8)(1.5,0){7}{\line(1,0){.7}}
\put(-4,-3){\tiny$y=-1/2$}\put(3,-11){\tiny$y=1/2-\ell$}\put(3,-16){\tiny $2\leq\ell\leq k$}\put(-4,21){\tiny$y=k-1/2$}\put(4,17){\tiny $y=k-3/2$}}
		%---5
		\put(100,15){\put(0,26){\tiny  Case (5)} \put(0,0){\line(1,0){4}}\put(4,0){\line(0,1){20}}\put(4,20){\line(-1,0){4}}\put(0,20){\line(0,-1){20}} \put(4,-4){\put(0,0){\line(1,0){4}}\put(4,0){\line(0,1){18}}\put(4,18){\line(-1,0){4}}\put(0,18){\line(0,-1){18}}} \multiput(4,20)(0.8,0){12}{\circle*{.01}}
\multiput(4,-8)(1.5,0){7}{\line(1,0){.7}}
\put(-4,-3){\tiny$y=-1/2$}\put(3,-11){\tiny$y=1/2-\ell$}
\put(3,-16){\tiny$\ell\geq k+1$}\put(-4,21){\tiny $y=k-1/2$}} }
	\end{picture}
	\caption{Decomposition of a polyomino in $\CP^u_k$, where
		the dotted (dashed) lines are the highest (lowest) position of
		top (bottom) cell of the second column} \label{figu1}
\end{figure}

Let $C^{u}_k(x,y,t) := \sum_{\nu\in \CP_k^{u}} x^{v(\nu)} y^{h(\nu)} t^{a(\nu)}.$
Define
$$C^{u}(x,y,t;z):=\sum_{k\geq1}C^{u}_k(x,y,t)z^{k-1} = \sum_{\nu\in \CP^{u}} x^{v(\nu)} y^{h(\nu)} t^{a(\nu)} z^{c(\nu)}.$$
Let $\nu\in\CP^u_k$. We point out that $\nu$ falls exclusively into one of the following cases (see Figure \ref{figu1}):
\begin{enumerate}
\item $\nu$ has one column;
\item $\nu$ has at least two columns such that
$b(\nu;1)<b(\nu;2)<u(\nu;2)\leq u(\nu;1)$; in this case $\nu$
without its first column is a nonempty polyomino in $\CP^{u}$;
\item $\nu$ has at least two columns such that $b(\nu;1)=b(\nu;2)$ and
$u(\nu;2)\leq u(\nu;1)$;
\item $\nu$ has at least two columns such that $b(\nu;2)<b(\nu;1)$,
$u(\nu;2)<u(\nu;1)$, and the second column has $\ell$ cells where
$2\leq \ell\leq k$;
\item $\nu$ has at least two columns such that $b(\nu;2)<b(\nu;1)$,
$u(\nu;2)\leq u(\nu;1)$, and the second column has $\ell\geq k+1$
cells.
\end{enumerate}

Thus, adding up the contributions of these cases, we have
\beq
C^{u}_k(x,y,t)&=&xy^kt^k+\sum_{l=1}^{k-1}(k-l)xy^{k-l}t^kC^{bu}_l(x,y,t)
+\sum_{l=1}^kxy^{k-l}t^kC^{u}_l(x,y,t)\\
&&+\sum_{l=2}^kx(y^{k-1}+\cdots+y^{k+1-l})t^kC^{u}_l(x,y,t)\\
&&+\sum_{l\geq k+1}x(y^{k-1}+\cdots+y+1)t^kC^{u}_l(x,y,t).
\feq
By multiplying by $z^{k-1}$, summing over $k\geq1$ and exchanging the order of the double sums, we have
\beq
C^{u}(x,y,t;z)&=&\frac{xyt}{1-ytz}+\sum_{l\geq1}\sum_{k\geq l+1}(k-l)xy^{k-l}t^kz^{k-1}C^{bu}_l(x,y,t)\\
&&+\sum_{l\geq1}\sum_{k\geq l}xy^{k-l}t^az^{k-1}C^{u}_l(x,y,t)\\
&&+\sum_{l\geq2}\sum_{k\geq l}x\frac{y^{k+1-l}-y^k}{1-y}t^kz^{k-1}C^{u}_l(x,y,t)\\
&&+\sum_{l\geq2}\sum_{k=1}^{l-1}x\frac{1-y^k}{1-y}t^kz^{k-1}C^{u}_l(x,y,t).
\feq
This results in the following equation
\beq
C^{u}(x,y,t;z)&=&\frac{xyt}{1-ytz}+\sum_{l\geq1}\frac{xyt^{l+1}z^l}{(1-ytz)^2}C^{bu}_l(x,y,t)
+\sum_{l\geq1}\frac{xt^bz^{b-1}}{1-tyz}C^{u}_l(x,y,t)\\
&&+\sum_{l\geq1}\frac{x(y-y^l)}{(1-y)(1-ytz)}t^lz^{l-1}C^{u}_l(x,y,t)\\
&&+\sum_{l\geq2}\frac{xt}{(1-tz)(1-ytz)}C^{u}_l(x,y,t)\\
&&-\sum_{l\geq2}\frac{x(y^ltz-ytz-y^l+1)}{(1-y)(1-tz)(1-ytz)}t^lz^{l-1}C^{u}_l(x,y,t).
\feq
Hence, the generating function $C^{u}(x,y,t;z)$ satisfies
\begin{align}\label{eqCPu1}
C^{u}(x,y,t;z)&=\frac{xyt}{1-ytz}+\frac{xyt^2z}{(1-ytz)^2}C^{bu}(x,y,t;tz)
+\frac{xt}{1-tyz}C^{u}(x,y,t;tz)\nonumber\\
&+\frac{xyt}{(1-y)(1-ytz)}C^{u}(x,y,t;tz)+\frac{xt}{(1-tz)(1-ytz)}C^{u}(x,y,t;1)\\
&-\frac{xt}{(1-y)(1-tz)}C^{u}(x,y,t;tz).\nonumber
\end{align}
Letting $t=1$ in \eqref{eqCPu1}, we obtain
\begin{align}\label{eqCPu2}
&\left(1+\frac{xz}{(1-z)(1-yz)}\right)C^{u}(x,y,1;z)\nonumber\\
&\qquad=\frac{xy}{1-yz}+\frac{xyz}{(1-yz)^2}C^{bu}(x,y,1;z)
+\frac{x}{(1-z)(1-yz)}C^{u}(x,y,1;1).
\end{align}
We solve this functional equation through an application of the kernel method (See \cite{Ban} for an introduction). To that end, let $z$ take the value
\beqn \label{z0_exm1}
z_0=z_0(x,y)=\frac{1+y-x-\sqrt{(1+y-x)^2-4y}}{2y},
\feqn
which is the root of $1+\frac{xz}{(1-z)(1-yz)}=0$. Then, by substitution in \eqref{eqCPu2}, we have
\beqn \label{onece_C1}
C^{u}(x,y,1;1)&=&y(z_0-1)+\frac{yz_0(z_0-1)}{1-yz_0}C^{bu}(x,y,1;z_0)\notag \\
&=&y(z_0-1)+\frac{xy^2z_0(z_0-1)}{(1-yz_0)^2-x},
\feqn
where we used Lemma \ref{lem1} for the last equality.
Hence, \eqref{onece_C1} along with \eqref{eqCPu2} and Lemma \ref{lem1} implies
\begin{lemma}\label{lem2}
The generating function $C^{u}(x,y,1;z)$ is given by
$$C^{u}(x,y,1;z)=\frac{xy(z_0-z)+\frac{x^2y^2z(1-z)}{(1-yz)^2-x}
-\frac{x^2y^2z_0(1-z_0)}{(1-yz_0)^2-x}}{(1-z)(1-yz)+xz},$$ where
$z_0$ is given by \eqref{z0_exm1}.
\end{lemma}
Note that $Z(x):=z_0(x,x)=\frac{1-\sqrt{1-4x}}{2x}$ is the generating function for the Catalan numbers $\frac{1}{n+1}\binom{2n}{n}$. The generating function for the number of polyominoes in $\CP^{u}$ with perimeter $2n$ is therefore given by
\beq
C^{u}(x,x,1;1)&=&x(Z(x)-1)+\frac{x^3Z(x)(Z(x)-1)}{(1-xZ(x))^2-x}\\
&=&\frac{x^2Z(x)}{2-Z(x)} =\frac{x^2}{\sqrt{1-4x}}=\sum_{n\geq2}\binom{2n-4}{n-2}x^n.
\feq
where we used the fact that $Z(x)=1+xZ^2(x)$. Consequently, the number of polyominoes in $\CP^{u}$ with perimeter $2n$ is given by $\binom{2n-4}{n-2}$.

\subsection{Enumeration in $\CP$}
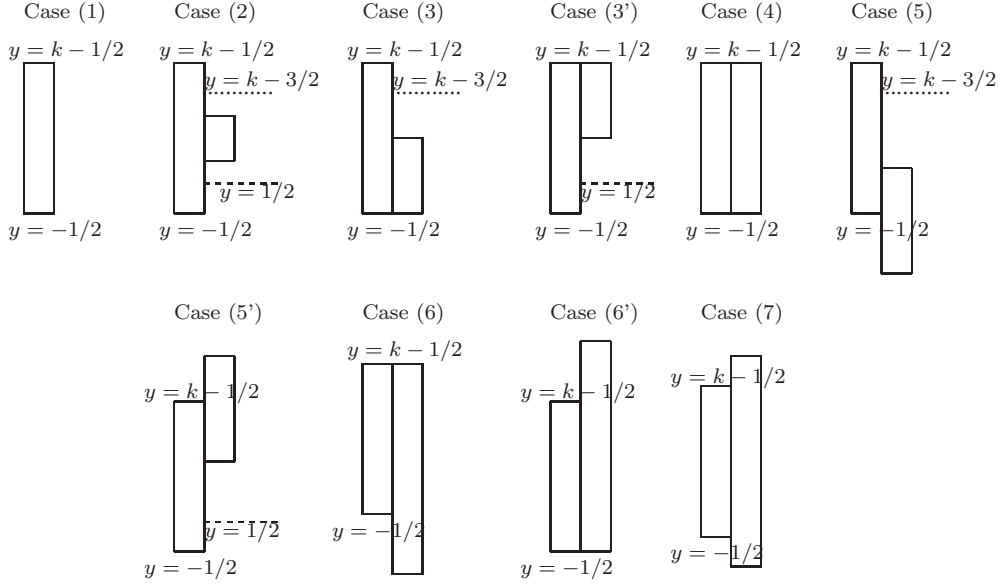
\begin{figure}[htp]
\begin{picture}(120,80)
	\setlength{\unitlength}{0.1cm} \put(0,-7){
\put(0,40){
		%---1
		\put(0,15){\put(0,26){\tiny Case (1)} \put(0,0){\line(1,0){4}}\put(4,0){\line(0,1){20}}\put(4,20){\line(-1,0){4}}\put(0,20){\line(0,-1){20}}
\put(-2,-3){\tiny$y=-1/2$}\put(-2,21){\tiny $y=k-1/2$}}
		%---2
		\put(20,15){\put(0,26){\tiny Case (2)} \put(0,0){\line(1,0){4}}\put(4,0){\line(0,1){20}}\put(4,20){\line(-1,0){4}}\put(0,20){\line(0,-1){20}}		\put(4,7){\put(0,0){\line(1,0){4}}\put(4,0){\line(0,1){6}}\put(4,6){\line(-1,0){4}}\put(0,6){\line(0,-1){6}}} \multiput(4,16)(0.8,0){12}{\circle*{.01}}
\multiput(4,4)(1.5,0){7}{\line(1,0){.7}}
\put(6,2){\tiny$y=1/2$}\put(-2,-3){\tiny$y=-1/2$}\put(-2,21){\tiny$y=k-1/2$}
\put(4,17){\tiny $y=k-3/2$}}
		%---3
		\put(45,15){\put(0,26){\tiny Case (3)} \put(0,0){\line(1,0){4}}\put(4,0){\line(0,1){20}}\put(4,20){\line(-1,0){4}}\put(0,20){\line(0,-1){20}} \put(4,0){\put(0,0){\line(1,0){4}}\put(4,0){\line(0,1){10}}\put(4,10){\line(-1,0){4}}\put(0,10){\line(0,-1){10}}}
\multiput(4,16)(0.8,0){12}{\circle*{.01}}
\put(-2,-3){\tiny$y=-1/2$}\put(-2,21){\tiny$y=k-1/2$}
\put(4,17){\tiny $y=k-3/2$}}
		%---31
		\put(70,15){\put(0,26){\tiny Case (3')} \put(0,0){\line(1,0){4}}\put(4,0){\line(0,1){20}}\put(4,20){\line(-1,0){4}}\put(0,20){\line(0,-1){20}}			\put(4,10){\put(0,0){\line(1,0){4}}\put(4,0){\line(0,1){10}}\put(4,10){\line(-1,0){4}}\put(0,10){\line(0,-1){10}}}
\multiput(4,4)(1.5,0){7}{\line(1,0){.7}}
\put(-2,-3){\tiny$y=-1/2$}\put(-2,21){\tiny$y=k-1/2$}
\put(4,2){\tiny$y=1/2$}}
		%---4
		\put(90,15){\put(0,26){\tiny Case (4)} \put(0,0){\line(1,0){4}}\put(4,0){\line(0,1){20}}\put(4,20){\line(-1,0){4}}\put(0,20){\line(0,-1){20}} \put(4,0){\put(0,0){\line(1,0){4}}\put(4,0){\line(0,1){20}}\put(4,20){\line(-1,0){4}}\put(0,20){\line(0,-1){20}}}
\put(-2,-3){\tiny$y=-1/2$}\put(-2,21){\tiny$y=k-1/2$}}
		%---5
		\put(110,15){\put(0,26){\tiny Case (5)} \put(0,0){\line(1,0){4}}\put(4,0){\line(0,1){20}}\put(4,20){\line(-1,0){4}}\put(0,20){\line(0,-1){20}}			\put(4,-8){\put(0,0){\line(1,0){4}}\put(4,0){\line(0,1){14}}\put(4,14){\line(-1,0){4}}\put(0,14){\line(0,-1){14}}}
\multiput(4,16)(0.8,0){12}{\circle*{.01}}
\put(-2,-3){\tiny$y=-1/2$}\put(-2,21){\tiny$y=k-1/2$}
\put(4,17){\tiny $y=k-3/2$}}
}
\put(-110,0){
		%---51
		\put(130,10){\put(0,31){\tiny Case (5')}	\put(0,0){\line(1,0){4}}\put(4,0){\line(0,1){20}}\put(4,20){\line(-1,0){4}}\put(0,20){\line(0,-1){20}} \put(4,12){\put(0,0){\line(1,0){4}}\put(4,0){\line(0,1){14}}\put(4,14){\line(-1,0){4}}\put(0,14){\line(0,-1){14}}}
\multiput(4,4)(1.5,0){7}{\line(1,0){.7}}
\put(-4,-3){\tiny$y=-1/2$}\put(-4,20.5){\tiny$y=k-1/2$}
\put(4,2){\tiny$y=1/2$}}
		%---6
		\put(155,15){\put(0,26){\tiny Case (6)} \put(0,0){\line(1,0){4}}\put(4,0){\line(0,1){20}}\put(4,20){\line(-1,0){4}}\put(0,20){\line(0,-1){20}}			\put(4,-8){\put(0,0){\line(1,0){4}}\put(4,0){\line(0,1){28}}\put(4,28){\line(-1,0){4}}\put(0,28){\line(0,-1){28}}}
\put(-4,-3){\tiny$y=-1/2$}\put(-2,21){\tiny$y=k-1/2$}}
		%---61
		\put(180,10){\put(0,31){\tiny Case (6')} \put(0,0){\line(1,0){4}}\put(4,0){\line(0,1){20}}\put(4,20){\line(-1,0){4}}\put(0,20){\line(0,-1){20}} \put(4,0){\put(0,0){\line(1,0){4}}\put(4,0){\line(0,1){28}}\put(4,28){\line(-1,0){4}}\put(0,28){\line(0,-1){28}}}
\put(-4,-3){\tiny$y=-1/2$}\put(-4,20.5){\tiny$y=k-1/2$}}
		%---7
		\put(200,12){\put(0,29){\tiny Case (7)} \put(0,0){\line(1,0){4}}\put(4,0){\line(0,1){20}}\put(4,20){\line(-1,0){4}}\put(0,20){\line(0,-1){20}}			\put(4,-4){\put(0,0){\line(1,0){4}}\put(4,0){\line(0,1){28}}\put(4,28){\line(-1,0){4}}\put(0,28){\line(0,-1){28}}}
\put(-4,-3){\tiny$y=-1/2$}\put(-4,20.5){\tiny$y=k-1/2$}} }}
	\end{picture}
	\caption{Decomposition of a nonempty polyomino in $\CP$, where the
		dotted (dashed) lines are the highest (lowest) position of top
		(bottom) cell of the second column}\label{figcp1}
\end{figure}

Let $C_k(x,y,t) := \sum_{\nu\in \CP_k} x^{v(\nu)} y^{h(\nu)} t^{a(\nu)}.$ Recall
\beq
C(x,y,t;z):=\sum_{k\geq1}C_k(x,y,t)z^{k-1}=\sum_{\nu\in \CP} x^{v(\nu)} y^{h(\nu)} t^{a(\nu)} z^{c(\nu)}.
\feq
Let $\nu\in\CP_k$. Then $\nu$ has either one column or at least two columns. In the later case, suppose that the number of cells in the second column is
$l$. With this convention, $\nu$ falls into one of the following cases. See Figure
\ref{figcp1} for a pictorial description of these cases.

\begin{itemize}
\item[(1)] $\nu$ has one column;
\item[(2)] $\nu$ has two columns such that $b(\nu;1)<b(\nu;2)<
u(\nu;2)< u(\nu;1)$; in this case $\nu$, without its first column,
is a nonempty polyomino in $\CP$;
\item[(3)] $\nu$ has two columns such that $b(\nu;1)=b(\nu;2)<
u(\nu;2)<u(\nu;1)$ (respectively, $u(\nu;1)=u(\nu;2)>
b(\nu;2)>b(\nu;1)$); in this case $\nu$, without its first column,
is a nonempty polyomino in $\CP^{u}$ (respectively, $\CP^{b}$);
\item[(4)] $\nu$ has two columns such that $b(\nu;2)=b(\nu;1)$ and
$u(\nu;2)=u(\nu;1)$; in this case $\nu$, without its first column,
is a nonempty polyomino in $\CP$;
\item[(5)] $\nu$ has two columns such that either $b(\nu;2)<b(\nu;1)$ and
$u(\nu;2)< u(\nu;1)$ or $b(\nu;2)>b(\nu;1)$ and
$u(\nu;2)> u(\nu;1)$, distinguishing two scenarios $l\geq
k+1$ and $2\leq l\leq k$; in this case $\nu$, without its first
column, is a nonempty polyomino in $\CP^u$ or $\CP^b$;
\item[(6)] $\nu$ has two columns such that either $b(\nu;2)<b(\nu;1)$ and
$u(\nu;2)=u(\nu;1)$ or $b(\nu;2)=b(\nu;1)$ and $u(\nu;2)>
u(\nu;1)$; in this case $\nu$, without its first column, is a
nonempty polyomino in $\CP$;
\item[(7)] $\nu$ has two columns such that $b(\nu;2)>b(\nu;1)$ and
$u(\nu;2)> u(\nu;1)$; in this case $\nu$, without its first column,
is a nonempty polyomino in $\CP$.
\end{itemize}
Considering these cases, we can naturally rewrite $C_k(x,y,t)$ as the sum of the following terms with the labels connecting them back to the corresponding cases; that is
\beq
&&C_k(x,y,t)=\underbrace{xy^kt^k}_{(1)}+\underbrace{\sum_{l=1}^{k-1}(k-1-l)xy^{k-l}t^kC_l^{bu}(x,y,t)}_{(2)}
+\underbrace{2\sum_{l=1}^{k-1}xy^{k-l}t^kC_l^{u}(x,y,t)}_{(3)}\\
&&\qquad \qquad  +\underbrace{xt^kC_k(x,y,t)}_{(4)}+\underbrace{2\sum_{l=2}^kx(y^{k-1}+\cdots+y^{k+1-l})t^kC_l^u(x,y,t)}_{(5)}\\
&&\qquad \qquad +\underbrace{2\sum_{l\geq
k+1}x(y^{k-1}+\cdots+y)t^kC_l^u(x,y,t)}_{(5)}
+\underbrace{2\sum_{l\geq k+1}xt^kC_l(x,y,t)}_{(6)}\\
&&\qquad \qquad +\underbrace{\sum_{l\geq k+1}(l-1-k)xt^kC_l(x,y,t)}_{(7)},
\feq
where we used the fact that by an easy symmetrical argument $C_l^{u}(x,y,t)=C_l^{b}(x,y,t)$ for any $l\geq1$. Multiplying $C_k(x,y,t)$ by $z^{k-1}$ and summing up over $k\geq1$, we obtain
\beq
&& C(x,y,t;z) =\\
&&\frac{xyt}{1-ytz}+x\sum_{k\geq1}\sum_{l=1}^{k-1}(k-1-l)y^{k-l}t^kz^{k-1}C_l^{bu}(x,y,t)
+2x\sum_{k\geq1}\sum_{l=1}^{k-1}y^{k-l}t^kz^{k-1}C_l^{u}(x,y,t)\\
&&+x\sum_{k\geq1}t^kz^{k-1}C_k(x,y,t)+2x\sum_{k\geq1}\sum_{l=2}^k\frac{y^{k+1-l}-y^k}{1-y}t^kz^{k-1}C_l^u(x,y,t)\\
&&+2x\sum_{k\geq1}\sum_{l\geq
k+1}\frac{y-y^k}{1-y}t^kz^{k-1}C_l^u(x,y,t)
+2x\sum_{k\geq1}\sum_{l\geq k+1}t^kz^{k-1}C_l(x,y,t)\\
&&+x\sum_{k\geq1}\sum_{l\geq k+1}(l-1-k)t^kz^{k-1}C_l(x,y,t).
\feq
Exchanging the order of double sums in all the terms gives
\beq
C(x,y,t;z)&=&\frac{xyt}{1-ytz}+\frac{xy^2t^3z^2}{(1-ytz)^2}C^{bu}(x,y,t;tz)
+\frac{2xyt^2z}{1-ytz}C^{u}(x,y,t;tz)\\
&&+xtC(x,y,t;tz)+\frac{2xyt}{(1-y)(1-ytz)}(C^u(x,y,t;tz)-C^u(x,y,t;ytz))\\
&&+\frac{2xyt}{(1-y)(1-tz)}(C^u(x,y,t;1)-C^u(x,y,t;tz))\\
&&-\frac{2xyt}{(1-y)(1-ytz)}(C^u(x,y,t;1)-C^u(x,y,t;ytz))\\
&&+\frac{2xt}{1-tz}(C(x,y,t;1)-C(x,y,t;tz))
+\frac{xt}{(1-tz)^2}C(x,y,t;tz)\\
&&-\frac{xt}{(1-tz)^2}C(x,y,t;1)
+\frac{xt}{1-tz}\frac{\partial}{\partial z}C(x,y,t;z)\mid_{z=1}.
\feq
This is equivalent to
\beqn\label{eqacp}
C(x,y,t;z)&=&\frac{xyt}{1-ytz}+\frac{xt^3z^2}{(1-tz)^2}C(x,y,t;tz)
+\frac{xt(1-2tz)}{(1-tz)^2}C(x,y,t;1)\nonumber\\
&&+\frac{xt}{1-tz}\frac{\partial}{\partial z}C(x,y,t;z)\mid_{z=1}
+\frac{xy^2t^3z^2}{(1-ytz)^2}C^{bu}(x,y,t;tz)\\
&&-\frac{2xyt^2z}{(1-ytz)(1-tz)}(tzC^{u}(x,y,t;tz)-\nonumber
C^u(x,y,t;1)).
\feqn
Setting $t=1$, we have
\beqn\label{eqacp1}
\left(1-\frac{xz^2}{(1-z)^2}\right)C(x,y,1;z)&=&\frac{xy}{1-yz}
+\frac{x(1-2z)}{(1-z)^2}C(x,y,1;1)\nonumber\\
&&+\frac{x}{1-z}\frac{\partial}{\partial z}C(x,y,1;z)\mid_{z=1}
+\frac{xy^2z^2}{(1-yz)^2}C^{bu}(x,y,1;z)\\
&&-\frac{2xyz}{(1-yz)(1-z)}(zC^{u}(x,y,1;z)-C^u(x,y,1;1)).\nonumber
\feqn
We apply the kernel method as we did in the previous section. To that end, we let $z=z_\pm(x,y)=\frac{1}{1\pm\sqrt{x}}$. By substituting these two values back in \eqref{eqacp1} we get the following two equations:
\beqn\label{eqacp111}
&&\frac{y(1-z_\pm)}{1-yz_\pm}+\frac{(1-2z_\pm)}{1-z_\pm}C(x,y,1;1)+\frac{\partial}{\partial z}C(x,y,1;z)\mid_{z=1}\nonumber\\
&&+\frac{y^2z_\pm^2(1-z_\pm)}{(1-yz_\pm)^2}C^{bu}(x,y,1;z_\pm)
-\frac{2yz_\pm}{1-yz_\pm}(z_\pm C^{u}(x,y,1;z_\pm)-C^u(x,y,1;1))=0.
\feqn
Subtracting these equations from each other, we obtain the following result

\begin{theorem}\label{thcp1}
We have
\beq
&&C(x,y,1;1)=\frac{2xy}{x-(1-y)^2}C^u(x,y,1;1)\\
&&-\frac{\sqrt{x}y}{(1+\sqrt{x})(1+\sqrt{x}-y)}C^u(x,y,1;z_+)
+\frac{\sqrt{x}y}{(1-\sqrt{x})(1-\sqrt{x}-y)}C^u(x,y,1;z_-)\\
&&+\frac{xy^2}{2(1+\sqrt{x})(1+\sqrt{x}-y)^2}C^{bu}(x,y,1;z_+)
+\frac{xy^2}{2(1-\sqrt{x})(1-\sqrt{x}-y)^2}C^{bu}(x,y,1;z_-)\\
&&-\frac{x(1-y)y}{x-(1-y)^2},
\feq
where $C^{bu}(x,y,1;z)$ and $C^u(x,y,1;z)$ are given by Lemma \ref{lem1} and Lemma \ref{lem2}, respectively.
\end{theorem}

Finally, setting $y=x$ in $C(x,y,1;1)$ gives the generating function \eqref{cp_perim_gf} and \eqref{eqfcpf}.

\section{Interior vertices}
\label{internal_ver}

In this section, we enumerate the convex polyominoes with respect the number of interior vertices and prove \eqref{total_int_exp}. To that end, define $F^{bu}(x,y,q;z)=\sum_{k\geq1}F^{bu}_k(x,y,q)z^{k-1}$, where
\beq
F^{bu}_k(x,y,q) := \sum_{\nu\in \CP^{bu}_k}x^{v(\nu)}y^{h(\nu)}q^{int(\nu)}.
\feq
Let $F^{u}_k(x,y,q)$, $F^{u}(x,y,q;z)$, $F_k(x,y,q)$, and $F(x,y,q;z)$ be defined in a similar fashion over $\CP^u$ and $\CP$.
We follow the arguments given in Section \ref{CP_decom} closely, however, we omit the details for the sake of space. First, we obtain the generating function over $\CP^{ub},$ that is
\beqn
F^{bu}(x,y,q;z)&=&\frac{xy}{1-yz}+ \frac{x}{(1-yz)^2}F^{bu}(x,y,q;qz),\label{eqaF1}
\feqn
Setting $q=1$ and solving for $F(x,y,1;z)$ in \eqref{eqaF1}, we have
\beqn
F^{bu}(x,y,1;z)=\frac{xy(1-yz)}{(1-yz)^2-x}.\label{eqbF11}
\feqn
Additionally, we take the derivative of \eqref{eqaF1} at $q=1$
\beqn
\frac{\partial}{\partial q}F^{bu}(x,y,q;z)\mid_{q=1}=
\frac{x^2y^2z((1-yz)^2y^2+x)}{((1-yz)^2-x)^3}.\label{eqbF2}
\feqn
Next, we calculate the generating function over $\CP^u$:
\beqn
F^{u}(x,y,q;z)&=&\frac{xy}{1-yz}+\frac{xyz}{(1-yz)^2}F^{bu}(x,y,q;qz)
-\frac{qxz}{(1-yz)(1-qz)}F^{u}(x,y,q;qz)\nonumber \\
&&+\frac{x}{(1-qz)(1-yz)}F^{u}(x,y,q;1),\label{eqaF22}
\feqn
To obtain $F^{u}(x,y,1;z)$, we set $q=1$ in \eqref{eqaF22}, use \eqref{eqbF11}, and solve by kernel method:
\beqn
F^{u}(x,y,1;z)&=&
\frac{x^2y}{\sqrt{(1+x-y)^2-4x}} \notag \\
&& +\frac{xy(1-z)(1-yz)(1-x-yz)}{(yz^2-(1-x+y)z+1)((1-yz)^2-x)}.\label{eqbF3}
\feqn
Furthermore, taking the derivative of \eqref{eqaF22} at $q=1$, we obtain
\beqn
&&\left(1-\frac{xz}{(1-z)(1-yz)}\right)\frac{\partial}{\partial q}F^{u}(x,y,q;z)\mid_{q=1} \nonumber \\
&&=\frac{xyz}{(1-yz)^2}\left(\frac{\partial}{\partial q}F^{bu}(x,y,q;z)\mid_{q=1}+z\frac{\partial}{\partial z}F^{bu}(x,y,1;z)\right)\nonumber\\
&&+\frac{xz}{(1-z)^2(1-yz)}(F^u(x,y,1;1)-zF^u(x,y,1;z))
+\frac{x}{1-yz}\frac{\partial}{\partial q}F^{u}(x,y,q;1)\mid_{q=1}\nonumber\\
&&-\frac{xz}{(1-z)(1-yz)}F^u(x,y,1;z)-\frac{xz^2}{(1-z)(1-yz)}\frac{\partial}{\partial z}F^{u}(x,y,1;z).\label{eqbF3q}
\feqn
Set $z=z_0$ where $z_0$ is the root of $\frac{xz}{(1-z)(1-yz)}=1.$ Then \eqref{eqaF22}-\eqref{eqbF3} imply
\beqn\label{eqbF4}
\frac{\partial}{\partial q}F^{u}(x,y,q;1)\mid_{q=1}
&=&\frac{xy(x^2-2xy+y^2+2x+2y-3)}{2\sqrt{(1+x-y)^2-4x}^3}\nonumber\\
&&+\frac{xy((1+x+y)(x-y)^2+4xy-5x-5y+3)}{2((1+x-y)^2-4x)^2}.
\feqn
We remark that one can obtain $\frac{\partial}{\partial q}F^{u}(x,y,q;z)\mid_{q=1}$ by substituting \eqref{eqbF4} into \eqref{eqbF3q}. We choose however to not include the explicit form of the generating function due to its length.
The last step is to obtain the generating function over $\CP:$

\beqn
F(x,y,q;z)&=&\frac{xy}{1-yz}+\frac{xy^2z^2}{(1-yz)^2}F^{bu}(x,y,q;qz)
-\frac{2qxyz^2}{(1-yz)(1-qz)}F^{u}(x,y,q;qz)\nonumber\\
&&+\frac{2xyz}{(1-yz)(1-qz)}F^u(x,y,q;1)+\frac{x}{(1-qz)^2}F(x,y,q;1)\nonumber\\
&&+\frac{q^2xz^2}{(1-qz)^2}F(x,y,q;qz)+\frac{x}{1-qz}\frac{\partial}{\partial z}F(x,y,q;z)\mid_{z=1}.\label{eqaF3}
\feqn
Recall that by definition of $F$ and $C$, $F(x,y,1;1) = C(x,y,1;1).$ Therefore, by \eqref{eqacp111} and Theorem \ref{thcp1},
\beqn
 F(x,y,1;1) &=& \frac{xy(2x^2(2-x)-x^2(3-x)(1-y)-2x(2-x)(1-y)^2+(1+x)(1-y)^3)}{((1+x-y)^2-4x)^2}\nonumber\\
&& -\frac{4x^2y^2}{\sqrt{(1+x-y)^2-4x}^3},\label{eqbF5}
\feqn
and
\beqn
\frac{\partial}{\partial z}F(x,y,1;z)\mid_{z=1}&=&\frac{xy^2(x^3-2x^2y+xy^2+x^2+3y^2-5x-6y+3)}{((1+x-y)^2-4x)^2} \nonumber \\
&&-\frac{2xy^2(1-x-y)}{\sqrt{(1+x-y)^2-4x}^3}.\label{eqbF6}
\feqn
Taking the derivative of \eqref{eqaF3} at $q=1$, we also obtain
\beqn \label{kernel5}
&&\left(1-\frac{xz^2}{(1-z)^2}\right)\frac{\partial}{\partial q}F(x,y,q;z)\mid_{q=1}\notag \\
&&=\frac{xy^2z^2}{(1-yz)^2} \left(\frac{\partial}{\partial q}F^{bu}(x,y,q;z)\mid_{q=1} +z\frac{\partial}{\partial z}F^{bu}(x,y,1;z)\right)+\frac{2xyz}{1-yz}\frac{\partial}{\partial q}F^u( x,y,q;z)\mid_{q=1} \notag\\
&&+\frac{2xyz^2}{1-yz}\frac{\partial}{\partial z}F^u(x,y,1;z)+\frac{2xyz^2}{(1-z)^2(1-yz)}(F^u(x,y,1;1)-F^u(x,y,1;z)) \notag\\
&&+\frac{2xyz}{(1-z)(1-yz)}\frac{\partial}{\partial q}F^u(x,y,q;1)\mid_{q=1} +\frac{xz^3}{(1-z)^2}\frac{\partial}{\partial z}F(x,y,1;z)-\frac{2x^2z}{(1-z)^3}F(x,y,1;1) \notag\\
&&+\frac{2x^2z}{(1-z)^3}F(x,y,1;z)+\frac{x(1-2z)}{(1-z)^2}\frac{\partial}{\partial q}F(x,y,q;1)\mid_{q=1}
+\frac{xz}{(1-z)^2}\frac{\partial}{\partial z}F(x,y,1;z)\mid_{z=1}  \notag \\
&&\qquad \qquad +\frac{x}{1-z}\frac{\partial^2}{\partial q\partial z}F(x,y,q;z)\mid_{q=z=1}.
\feqn
Given \eqref{eqbF2}-\eqref{eqbF3}, \eqref{eqbF4}, \eqref{eqbF5}, and \eqref{eqbF6}, we solve the equation \eqref{kernel5} using the kernel method. By substituting $z=z_+(x,y)=\frac{1}{1+\sqrt{x}}$ and $z=z_-(x,y)=\frac{1}{1-\sqrt{x}},$ we have a system of two equations with two variables $\frac{\partial}{\partial q}F(x,y,q;1)\mid_{q=1}$ and $\frac{\partial^2}{\partial q\partial z}F(x,y,q;z)\mid_{q=z=1}$. Solving this system, we get
\begin{theorem} \label{internal_thm}
The generating function $\frac{\partial}{\partial q}F(x,y,q;1)\mid_{q=1}$ is given by
\beq
&&\frac{8x^3y^2(4x(x^2+6x-1)+4x(5x+11)(y-x-1)+(6x^2+29x+17)(y-x-1)^2)}{((1+x-y)^2-4x)^4}\\
&&+\frac{x^2y^2(y-x-1)^3(40x(2x+1)+2(x^2+11x-8)(y-x-1)+(y-x-1)^2(y+x-1))}{((1+x-y)^2-4x)^4}\\
&&+\frac{4x^2y^2(4-3(x+y)-(x-y)^2)}{\sqrt{(1+x-y)^2-4x}^5}.
\feq
\end{theorem}

Theorem~\ref{internal_thm} implies that that the generating function for the total number of interior vertices in all convex polyominoes with perimeter $2n$ is
$$\frac{\partial}{\partial q}F(x,x,q;1)\mid_{q=1}
=\frac{(32x^4+80x^3-230x^2+116x-15)x^4}{(1-4x)^4}-\frac{8(3x-2)x^4}{\sqrt{1-4x}^5}.$$
The formula \eqref{total_int_exp} is then obtained by extracting the coefficient of $x^{2n}$ in this generating function.

\section{Boundary vertices of certain degrees}
\label{boundary_ver}

In this section, we will count the number of polyominoes with respect to boundary vertices of various degrees and prove \eqref{degree2_average}, \eqref{degree3_asymp}, \eqref{bv1-asym}, and \eqref{bv2-asym}. To that end, define
\beq
D_k(x,y,q_2,q_3,q_4) := \sum_{\nu \in \CP_k} x^{v(\nu)}y^{h(\nu)}q_2^{d_2(\nu)}q_3^{d_3(\nu)}q_4^{d_4(\nu)}.
\feq
where $q_i$ marks the boundary vertices of degree $i$. Similarly, let
\beqn \label{Dxy-def}
D(x,y,q_2,q_3,q_4;z):=\sum_{k\geq1}D_k(x,y,q_2,q_3,q_4)z^{k-1} = \sum_{\nu \in \CP} x^{v(\nu)}y^{h(\nu)}q_2^{d_2(\nu)}q_3^{d_3(\nu)}q_4^{d_4(\nu)}z^{c(\nu)}.
\feqn
Recall that for each $\nu\in\CP$
\beq
d_4(\nu)&=&d_2(\nu)-4 \\
h(\nu) + v(\nu) &=& \frac{1}{2}(d_2(\nu)+d_3(\nu)+d_4(\nu)).
\feq
Thus, from \eqref{Dxy-def} we get
\beqn\label{eqaD1}
D(q_3^2x,q_3^2y,q_2q_4/q_3^2,1,1;z)=\frac{q_4^4}{q_3^4}D(x,y,q_2,q_3,q_4;z).
\feqn
Therefore, it is adequate to study the generating function $E(x,y,q;z):=\sum_{k\geq1}E_k(x,y,q)z^{k-1}$ where
\beq
E_k(x,y,q) := \sum_{\nu \in \CP_k} x^{v(\nu)}y^{h(\nu)}q^{d_2(\nu)}
\feq
Similarly, we define the generating functions $E^{bu}(x,y,q;z)$ and $E^u(x,y,q;z)$ over $\CP^{bu}$ and $\CP^u$.
Following similar steps as the one given in the previous sections, we have
\beqn
E^{bu}(x,y,q;z)&=&\frac{q^4xy}{1-yz}+\frac{q^2xy^2z^2}{(1-yz)^2}E^{bu}(x,y,q;z)\nonumber\\
&&+\frac{2qxyz}{1-yz}E^{bu}(x,y,q;z)+xE^{bu}(x,y,q;z),\label{eqbD1}\\
E^{u}(x,y,q;z)&=&\frac{q^4xy}{1-yz}+\frac{q^2xy^2z^2}{(1-yz)^2}E^{bu}(x,y,q;z)
+\frac{qxyz}{1-yz}E^{bu}(x,y,q;z)\nonumber\\
&&+\frac{qxyz}{1-yz}E^{u}(x,y,q;z)+xE^{u}(x,y,q;z)
+\frac{q^2xy}{(1-y)(1-yz)}E^u(x,y,q;z)\nonumber\\
&&+\frac{q^2xy}{(1-y)(1-z)}(E^u(x,y,q;1)-E^u(x,y,q;z))\nonumber\\
&&-\frac{q^2xy}{(1-y)(1-yz)}E^u(x,y,q;1)+\frac{qx}{1-z}(E^u(x,y,q;1)-E^u(x,y,q;z)),\label{eqbD2}\\
E(x,y,q;z)&=&\frac{q^4xy}{1-yz}+\frac{q^2xy^2z^2}{(1-yz)^2}E^{bu}(x,y,q;z)
+\frac{qxyz}{1-yz}E^{u}(x,y,q;z)\nonumber\\
&&+\frac{qxyz}{1-yz}E^{u}(x,y,q;z)+xE(x,y,q;z)
+\frac{2q^2xy}{(1-y)(1-yz)}E^u(x,y,q;z)\nonumber\\
&&+\frac{2q^2xy}{(1-y)(1-z)}(E^u(x,y,q;1)-E^u(x,y,q;z))\nonumber\\
&&-\frac{2q^2xy}{(1-y)(1-yz)}E^u(x,y,q;1)
+\frac{2qxz}{1-z}(E(x,y,q;1)-E(x,y,q;z))\nonumber\\
&&+\frac{q^2x}{(1-z)^2}(E(x,y,q;z)-E(x,y,q;1))+\frac{q^2x}{1-z}\frac{\partial}{\partial}E(x,y,q;z)\mid_{z=1}.\label{eqbD3}
\feqn
By \eqref{eqbD1}, we have
\beqn
E^{bu}(x,y,q;z)&=\frac{q^4xy(1-yz)}{(1-yz)^2-x(1+(q-1)yz)^2}.\label{eqbD4}
\feqn
Hence, \eqref{eqbD4} along with \eqref{eqbD2} imply
\beqn
&&\left(1-\frac{x(z+q-1)(1+(q-1)yz)}{(1-yz)(1-z)}\right)E^{u}(x,y,q;z)\nonumber\\
&&\qquad=\frac{qx(1+(q-1)yz)}{(1-yz)(1-z)}E^u(x,y,q;1)+\frac{q^4xy(1-x-yz-(q-1)xyz)}{(1-yz)^2-x(1+(q-1)yz)^2}.\label{eqbD5}
\feqn
Once again we apply the kernel method to solve this equation. To that end, we substitute
$$z=z_0=\frac{1+y-x-(q-1)^2xy-\sqrt{(1+y-x-(q-1)^2xy)^2-4y(1+(q-1)x)^2}}{2(1+(q-1)x)}$$
into \eqref{eqbD5}. This substitution implies
\beqn
E^u(x,y,q;1)=-\frac{q^3y(1-x-yz_0-(q-1)xyz_0)(1-yz_0)(1-z_0)}{((1-yz_0)^2-x(1+(q-1)yz_0)^2)(1+(q-1)yz_0)}.\label{eqbD6}
\feqn
Inserting \eqref{eqbD6} into \eqref{eqbD5}, and solving for $E^u(x,y,q;z)$, we obtain
\beqn\label{eqbD7}
&&E^{u}(x,y,q;z)=\frac{q^4xy(1-x-yz-(q-1)xyz)}
{(1-x\frac{(z+q-1)(1+(q-1)yz)}{(1-yz)(1-z)})((1-yz)^2-x(1+(q-1)yz)^2)}\nonumber\\
&&-\frac{q^4xy(1+(q-1)yz)(1-x-yz_0-(q-1)xyz_0)\frac{(1-yz_0)(1-z_0)}{(1-yz)(1-z)}}{(1-x\frac{(z+q-1)(1+(q-1)yz)}{(1-yz)(1-z)})((1-yz_0)^2-x(1+(q-1)yz_0)^2)(1+(q-1)yz_0)}.
\feqn
Finally by \eqref{eqbD3}, we have
\beqn
&&\left(1-\frac{x(q-1+z)^2}{(1-z)^2}\right)E(x,y,q;z)\nonumber\\
&&=\frac{q^4xy}{1-yz}+\frac{q^2xy^2z^2}{(1-yz)^2}E^{bu}(x,y,q;z)
-\frac{2qxyz(q-1+z)}{(1-yz)(1-z)}E^{u}(x,y,q;z)\nonumber\\
&&+\frac{2q^2xyz}{(1-z)(1-yz)}E^u(x,y,q;1)
-\frac{qx(q-2+2z)}{(1-z)^2}E(x,y,q;1)\nonumber\\
&&+\frac{q^2x}{1-z}\frac{\partial}{\partial}E(x,y,q;z)\mid_{z=1}.\label{eqbD8}
\feqn
We now substitute $z=z_+=1-\frac{q\sqrt{x}}{1+\sqrt{x}}$ and $z=z_-=1+\frac{q\sqrt{x}}{1-\sqrt{x}}$ into \eqref{eqbD8}, and obtain two equations with two variables $E(x,y,q;1)$ and $\frac{\partial}{\partial}E(x,y,q;z)\mid_{z=1}$. Note that the expressions of $E^{bu}(x,y,q;z_\pm)$, $E^u(x,y,q;1)$, $E^{u}(x,y,q;z_\pm)$ are given by \eqref{eqbD4}, \eqref{eqbD6} and \eqref{eqbD7}, respectively. Solving this system of equations results in Theorem \ref{thdeg2}.
\begin{theorem}\label{thdeg2}
The generating function $E(x,y,q;1)$ is
\beq
&&-\frac{q^4x^2y^2(2+(x+y+1)(q-1)-xy(q-1)^3)^2}{\sqrt{(1+y-x-(q-1)^2xy)^2-4y(1+(q-1)x)^2}^3}\\
&&+\frac{q^4xy(-x^3y+2x^2y^2-xy^3-x^3-x^2y-xy^2-y^3+3x^2+5xy+3y^2-3x-3y+1)}{((1+y-x-(q-1)^2xy)^2-4y(1+(q-1)x)^2)^2}\\
&&-\frac{8q^4x^2y^2(q-1)(x^2-xy+y^2-2x-2y+1)}{((1+y-x-(q-1)^2xy)^2-4y(1+(q-1)x)^2)^2}\\
&&-\frac{q^4x^2y^2(q-1)^2(x^3+x^2y+xy^2+y^3+x^2-44xy+y^2-5x-5y+3)}{((1+y-x-(q-1)^2xy)^2-4y(1+(q-1)x)^2)^2}\\
&&+\frac{q^4x^3y^3(q-1)^3(8(2y+2x+3)+(3x^2+5xy+3y^2+5x+5y+4)(q-1))}{((1+y-x-(q-1)^2xy)^2-4y(1+(q-1)x)^2)^2}\\
&&-\frac{q^4x^4y^4(q-1)^5(8+3(x+y+1)(q-1)-xy(q-1)^3)}{((1+y-x-(q-1)^2xy)^2-4y(1+(q-1)x)^2)^2}.
\feq
\end{theorem}

As a corollary to Theorem \ref{thdeg2}, we get
\beq
\frac{\partial}{\partial q}E(x,x,q;1)\mid_{q=1}
=\frac{4x^2(8x^5-26x^4+56x^3-37x^2+10x-1)}{(4x-1)^3}-\frac{4x^4(4x^2-18x+5)}{\sqrt{1-4x}^5}.
\feq
Thus, by finding the coefficient of $x^n$, we arrive at \eqref{degree2_total}. Dividing by \eqref{eqfcpf} we obtain \eqref{degree2_average}. More importantly, \eqref{eqaD1} implies that

\beqn\label{eqaD11a}
D(x,y,q_2,q_3,q_4;1)=\frac{q_3^4}{q_4^4}E(q_3^2x,q_3^2y,q_2q_4/q_3^2; 1),
\feqn
where $E(x,y,q;1)$ is given by Theorem \ref{thdeg2}. In particular, the generating function $D(1,1,q,q,1;1)$, enumerating convex polyominoes with respect to boundary vertices of degrees two and three, is given by
\begin{align} \label{deg3_gf}
&\frac{(q^{12}-6q^{11}+11q^{10}-6q^9+6q^8-12q^7-4q^6+2q^5+9q^4+6q^3-4q^2-2q+1)q^4}{(q^2+q+1)^2(q^2-3q+1)^2(q^2-q-1)^2} \nonumber \\
&+\frac{(q^3-2q^2-1)^2q^6}{(q^2-q-1)\sqrt{(q^2+q+1)(q^2-3q+1)}^3}.
\end{align}
Then, the asymptotic result \eqref{degree3_asymp} is obtained by examining the coefficient of $q^n$.  Note that by \eqref{eqaD11a}, we have that $$q^4D(1,1,1,q,q;1)=D(1,1,q,q,1;1).$$
Hence, the number of polyominoes $\nu$ in $\CP$ with $d_2(\nu)+d_3(\nu)=n$ is equal to the number of polyominoes $\nu$ in $\CP$ with $d_3(\nu)+d_4(\nu)=n-4$, for all $n\geq4$.

Moreover, by considering the generating function $D(1,1,q,q,p;1)$, we obtain that the total number of all boundary vertices of degrees four over all polyominoes $\nu\in \CP$ with $d_2(\nu)+d_3(\nu)=n$ is asymptotically given by \eqref{bv1-asym}. This, in addition, implies that the expected number of the boundary vertices of degree four over all polyominoes in$\CP$ with  $d_2(\pi)+d_3(\pi)=n$ is asymptotically \eqref{bv2-asym}.

\section{Outer site perimeter of convex polyominoes} \label{outer_section}
In this section, we study the outer-site perimeter of convex polyominoes and prove  \eqref{tot_out_resp_perimeter}, \eqref{avg_out_resp_perimeter}, \eqref{avg_perim_res_outer} and \eqref{outer_asympt}.  To that end, let $J^{bu}_k(x,y,q) := \sum_{\nu\in \CP_k^{bu}} x^{v(\nu)} y^{h(\nu)} q^{o(\nu)}.$
Marking the cells in the first column with $z$ and summing over $k,$ we then define
$$J^{bu}(x,y,q;z):=\sum_{k\geq1}J^{bu}_k(x,y,q)z^{k-1} = \sum_{\nu\in \CP^{bu}} x^{v(\nu)} y^{h(\nu)} q^{o(\nu)} z^{c(\nu)}.$$ Similarly, we define the generating functions $J_k^{u}(x,y,q)$, $J^{u}(x,y,q;z)$, $J_k(x,y,q)$, and $J(x,y,q;z)$ over $\CP^u_k,$ $\CP^u$, $\CP_k,$ and $\CP$, respectively. \\

\subsection{Enumeration over $\CP^{ub}$}

Let $\nu\in\CP^{bu}_k$. Then, $\nu$ falls exclusively into one of the following cases (see Figure \ref{figd1}):
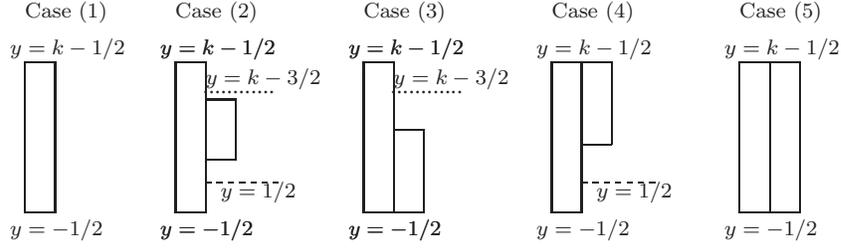
\begin{figure}[htp]
	\begin{picture}(100,30)
	\setlength{\unitlength}{0.1cm}\put(0,-13){
		%---1
		\put(0,15){
			\put(0,26){\tiny  Case (1)} \put(0,0){\line(1,0){4}}\put(4,0){\line(0,1){20}}\put(4,20){\line(-1,0){4}}\put(0,20){\line(0,-1){20}}
			\put(-2,-3){\tiny$y=-1/2$}\put(-2,21){\tiny$y=k-1/2$}}
		%---2
		\put(20,15){
			\put(0,26){\tiny  Case (2)} \put(0,0){\line(1,0){4}}\put(4,0){\line(0,1){20}}\put(4,20){\line(-1,0){4}}\put(0,20){\line(0,-1){20}}
			\put(4,7){\put(0,0){\line(1,0){4}}\put(4,0){\line(0,1){8}}\put(4,8){\line(-1,0){4}}\put(0,8){\line(0,-1){8}}}           \put(-2,-3){\tiny$y=-1/2$}\put(-2,21){\tiny$y=k-1/2$}
			\multiput(4,16)(0.8,0){12}{\circle*{.01}}
			\multiput(4,4)(1.5,0){7}{\line(1,0){.7}}
			\put(6,2){\tiny$y=1/2$}\put(-2,-3){\tiny$y=-1/2$}\put(-2,21){\tiny$y=k-1/2$}
			\put(4,17){\tiny$y=k-3/2$}}
		%---3
		\put(45,15){
			\put(0,26){\tiny  Case (3)} \put(0,0){\line(1,0){4}}\put(4,0){\line(0,1){20}}\put(4,20){\line(-1,0){4}}\put(0,20){\line(0,-1){20}}
			\put(4,0){\put(0,0){\line(1,0){4}}\put(4,0){\line(0,1){11}}\put(4,11){\line(-1,0){4}}\put(0,11){\line(0,-1){11}}}           \put(-2,-3){\tiny$y=-1/2$}\put(-2,21){\tiny$y=k-1/2$}
			\multiput(4,16)(0.8,0){12}{\circle*{.01}}
			\put(-2,-3){\tiny$y=-1/2$}\put(-2,21){\tiny$y=k-1/2$}
			\put(4,17){\tiny$y=k-3/2$}}
		%---4
		\put(70,15){
			\put(0,26){\tiny  Case (4)} \put(0,0){\line(1,0){4}}\put(4,0){\line(0,1){20}}\put(4,20){\line(-1,0){4}}\put(0,20){\line(0,-1){20}}
			\put(4,9){\put(0,0){\line(1,0){4}}\put(4,0){\line(0,1){11}}\put(4,11){\line(-1,0){4}}\put(0,11){\line(0,-1){11}}}           \put(-2,-3){\tiny$y=-1/2$}\put(-2,21){\tiny$y=k-1/2$}
			\multiput(4,4)(1.5,0){7}{\line(1,0){.7}}
			\put(6,2){\tiny$y=1/2$}}
		%---5
		\put(95,15){
			\put(0,26){\tiny  Case (5)} \put(0,0){\line(1,0){4}}\put(4,0){\line(0,1){20}}\put(4,20){\line(-1,0){4}}\put(0,20){\line(0,-1){20}}
			\put(4,0){\put(0,0){\line(1,0){4}}\put(4,0){\line(0,1){20}}\put(4,20){\line(-1,0){4}}\put(0,20){\line(0,-1){20}}}           \put(-2,-3){\tiny$y=-1/2$}\put(-2,21){\tiny$y=k-1/2$}}
	}
	\end{picture}
	\caption{Decomposition of a polyomino in $\CP_k^{bu}$}\label{figd1}
\end{figure}
\begin{itemize}
	\item[Case (1)] It has one column;
	\item[Case (2)] It has at least two columns such that $u(\nu,2)<u(\nu,1)$ and $b(\nu,2)>b(\nu,1)$;
	\item[Case (3)] It has at least two columns such that $u(\nu,2)<u(\nu,1)$ and $b(\nu,2)=b(\nu,1)$;
	\item[Case (4)] It has at least two columns such that $u(\nu,2)=u(\nu,1)$ and $b(\nu,2)>b(\nu,1)$;
	\item[Case (5)] It has at least two columns such that $u(\nu,2)=u(\nu,1)$ and $b(\nu,2)=b(\nu,1)$.
\end{itemize}
We then write a recurrence for $J^{bu}_k(x,y,q)$, by adding up the contribution of each of the above cases; that is
\beq
J^{bu}_k(x,y,q)&=&\underbrace{q^{2k+2}xy^k}_{(1)}
+\underbrace{\sum_{l=1}^{k-1}(k-1-l)xy^{k-l}q^{2k-2l}J^{bu}_l(x,y,q)}_{(2)}\\
&&+\underbrace{2\sum_{l=1}^{k-1}xy^{k-l}q^{2k-2l+1}J^{bu}_l(x,y,q)}_{(3)+(4)}
+\underbrace{q^2xJ_k(x,y,q)}_{(5)}.
\feq
Multiplying by $z^{k-1}$ and summing over $k\geq1$, we obtain
\beq
J^{bu}(x,y,q;z)&=&\frac{q^4xy}{1-q^2yz}+\frac{q^4xy^2z^2}{(1-q^2yz)^2}J^{bu}(x,y,q;z)\\
&&+\frac{2q^3xyz}{1-q^2yz}J^{bu}(x,y,q;z)+q^2xJ^{bu}(x,y,q;z),
\feq
We solve this equation for $J^{bu}(x,y,q;z)$ and get

\begin{lemma}\label{lem1_outer}
	The generating function $J^{bu}(x,y,q;z)$ is given by
	$$\frac{q^4xy(1-q^2yz)}{(1-q^2yz)^2-q^2(1+qyz-q^2yz)^2x}.$$
\end{lemma}

\subsection{Enumeration over $\CP^{u}$}

Let $\nu\in\CP^{u}_k$. We point out that $\nu$ falls exclusively into one of the following cases (see Figure \ref{figd2}):
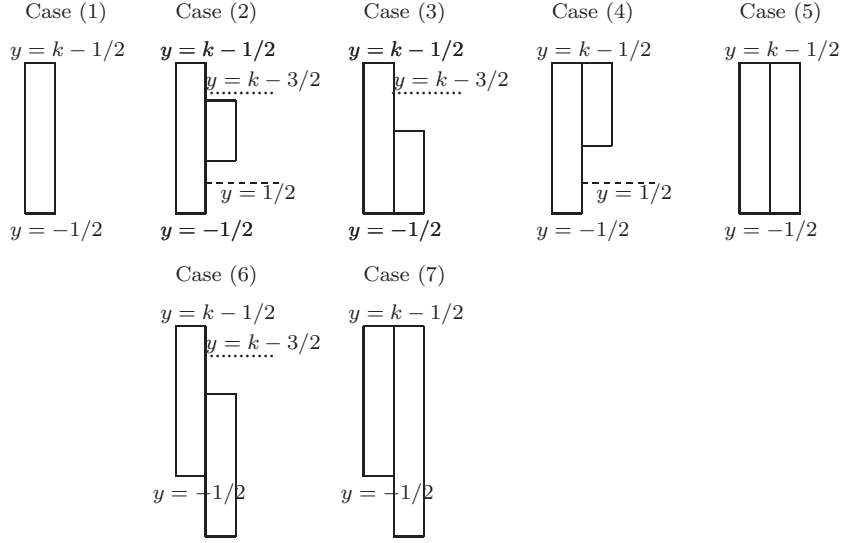
\begin{figure}[htp]
	\begin{picture}(100,73)
	\setlength{\unitlength}{0.1cm}\put(0,30){
		%---1
		\put(0,15){
			\put(0,26){\tiny  Case (1)} \put(0,0){\line(1,0){4}}\put(4,0){\line(0,1){20}}\put(4,20){\line(-1,0){4}}\put(0,20){\line(0,-1){20}}
			\put(-2,-3){\tiny$y=-1/2$}\put(-2,21){\tiny$y=k-1/2$}}
		%---2
		\put(20,15){
			\put(0,26){\tiny  Case (2)} \put(0,0){\line(1,0){4}}\put(4,0){\line(0,1){20}}\put(4,20){\line(-1,0){4}}\put(0,20){\line(0,-1){20}}
			\put(4,7){\put(0,0){\line(1,0){4}}\put(4,0){\line(0,1){8}}\put(4,8){\line(-1,0){4}}\put(0,8){\line(0,-1){8}}}           \put(-2,-3){\tiny$y=-1/2$}\put(-2,21){\tiny$y=k-1/2$}
			\multiput(4,16)(0.8,0){12}{\circle*{.01}}
			\multiput(4,4)(1.5,0){7}{\line(1,0){.7}}
			\put(6,2){\tiny$y=1/2$}\put(-2,-3){\tiny$y=-1/2$}\put(-2,21){\tiny$y=k-1/2$}
			\put(4,17){\tiny$y=k-3/2$}}
		%---3
		\put(45,15){
			\put(0,26){\tiny  Case (3)} \put(0,0){\line(1,0){4}}\put(4,0){\line(0,1){20}}\put(4,20){\line(-1,0){4}}\put(0,20){\line(0,-1){20}}
			\put(4,0){\put(0,0){\line(1,0){4}}\put(4,0){\line(0,1){11}}\put(4,11){\line(-1,0){4}}\put(0,11){\line(0,-1){11}}}           \put(-2,-3){\tiny$y=-1/2$}\put(-2,21){\tiny$y=k-1/2$}
			\multiput(4,16)(0.8,0){12}{\circle*{.01}}
			\put(-2,-3){\tiny$y=-1/2$}\put(-2,21){\tiny$y=k-1/2$}
			\put(4,17){\tiny$y=k-3/2$}}
		%---4
		\put(70,15){
			\put(0,26){\tiny  Case (4)} \put(0,0){\line(1,0){4}}\put(4,0){\line(0,1){20}}\put(4,20){\line(-1,0){4}}\put(0,20){\line(0,-1){20}}
			\put(4,9){\put(0,0){\line(1,0){4}}\put(4,0){\line(0,1){11}}\put(4,11){\line(-1,0){4}}\put(0,11){\line(0,-1){11}}}           \put(-2,-3){\tiny$y=-1/2$}\put(-2,21){\tiny$y=k-1/2$}
			\multiput(4,4)(1.5,0){7}{\line(1,0){.7}}
			\put(6,2){\tiny$y=1/2$}}
		%---5
		\put(95,15){
			\put(0,26){\tiny  Case (5)} \put(0,0){\line(1,0){4}}\put(4,0){\line(0,1){20}}\put(4,20){\line(-1,0){4}}\put(0,20){\line(0,-1){20}}
			\put(4,0){\put(0,0){\line(1,0){4}}\put(4,0){\line(0,1){20}}\put(4,20){\line(-1,0){4}}\put(0,20){\line(0,-1){20}}}           \put(-2,-3){\tiny$y=-1/2$}\put(-2,21){\tiny$y=k-1/2$}}
	}
	%----6
	\put(20,10){
		\put(0,26){\tiny  Case (6)} \put(0,0){\line(1,0){4}}\put(4,0){\line(0,1){20}}\put(4,20){\line(-1,0){4}}\put(0,20){\line(0,-1){20}}
		\put(4,-8){\put(0,0){\line(1,0){4}}\put(4,0){\line(0,1){19}}\put(4,19){\line(-1,0){4}}\put(0,19){\line(0,-1){19}}}
		\multiput(4,16)(0.8,0){12}{\circle*{.01}}
		\put(-2,21){\tiny$y=k-1/2$}
		\put(-3,-3){\tiny$y=-1/2$}\put(4,17){\tiny$y=k-3/2$}}
	%----7
	\put(45,10){
		\put(0,26){\tiny  Case (7)} \put(0,0){\line(1,0){4}}\put(4,0){\line(0,1){20}}\put(4,20){\line(-1,0){4}}\put(0,20){\line(0,-1){20}}
		\put(4,-8){\put(0,0){\line(1,0){4}}\put(4,0){\line(0,1){28}}\put(4,28){\line(-1,0){4}}\put(0,28){\line(0,-1){28}}}
		\put(-3,-3){\tiny$y=-1/2$}\put(-2,21){\tiny$y=k-1/2$}}
	\end{picture}
	\caption{Decomposition of a polyomino in $\CP_k^{u}$}\label{figd2}
\end{figure}
\begin{itemize}
	\item[Case (1)] It has one column;
	\item[Case (2)] It has at least two columns such that $u(\nu,2)<u(\nu,1)$ and $b(\nu,2)>b(\nu,1)$;
	\item[Case (3)] It has at least two columns such that $u(\nu,2)<u(\nu,1)$ and $b(\nu,2)=b(\nu,1)$;
	\item[Case (4)] It has at least two columns such that $u(\nu,2)=u(\nu,1)$ and $b(\nu,2)>b(\nu,1)$;
	\item[Case (5)] It has at least two columns such that $u(\nu,2)=u(\nu,1)$ and $b(\nu,2)=b(\nu,1)$;
	\item[Case (6)] It has at least two columns such that $u(\nu,2)<u(\nu,1)$ and $b(\nu,2)<b(\nu,1)$;
	\item[Case (7)] It has at least two columns such that $u(\nu,2)=u(\nu,1)$ and $b(\nu,2)<b(\nu,1)$.
\end{itemize}

Considering the contribution of each of these cases, we then write $J^{u}_k(x,y,q)$ recursively; that is,
\beq
J^{u}_k(x,y,q)&=&\underbrace{q^{2k+2}xy^k}_{(1)}
+\underbrace{\sum_{l=1}^{k-1}(k-1-l)xy^{k-l}q^{2k-2l}J^{bu}_l(x,y,q)}_{(2)}\\
&&+\underbrace{\sum_{l=1}^{k-1}xy^{k-l}q^{2k-2l+1}(J^{bu}_l(x,y,q)+J^u_l(x,y,q))}_{(3)+(4)}
+\underbrace{q^2xJ_k^u(x,y,q)}_{(5)}\\
&&+\underbrace{\sum_{l=2}^k\frac{xy^{k-l+1}q^{2k-2l+2}-xy^kq^{2k}}{1-q^2y}J_l^u(x,y,q)
	+\sum_{l\geq k+1}\frac{xyq^2-xy^kq^{2k}}{1-q^2y}J_l^u(x,y,q)}_{(6)}\\
&&+\underbrace{\sum_{l\geq k+1}xqJ_l^u(x,y,q)}_{(7)}.
\feq
We multiply by $z^{k-1}$ and sum over $k\geq1$. Exchanging the order of each double sum, we obtain
\beq
J^{u}(x,y,q;z)&=&\frac{q^4xy}{1-q^2yz}+\frac{q^4xy^2z^2}{(1-q^2yz)^2}J^{bu}(x,y,q;z)\\
&&+\frac{q^3xyz}{1-q^2yz}(J^{bu}(x,y,q;z)+J^u(x,y,q;z))+q^2xJ^{u}(x,y,q;z)\\
&&-\frac{qx(1+q(1-q)yz)}{(1-z)(1-q^2yz)}(J^u(x,y,q;z)-J^u(x,y,q;1)).
\feq
This is equivalent to
\begin{align}\label{eqCu1}
&\left(1-\frac{qx(1-q+qz)(1+q(1-q)yz)}{(1-z)(1-q^2yz)}\right)J^{u}(x,y,q;z)\nonumber\\
&\qquad=\frac{q^4xy}{1-q^2yz}+\frac{q^4xy^2z^2}{(1-q^2yz)^2}J^{bu}(x,y,q;z)
+\frac{q^3xyz}{1-q^2yz}J^{bu}(x,y,q;z)\\
&\qquad+\frac{qx(1+q(1-q)yz)}{(1-z)(1-q^2yz)}J^u(x,y,q;1).\nonumber
\end{align}
To solve this functional equation, we apply the kernel method. To that end, let $z$ take the value
\begin{align} \label{z0ker1}
z_0=\frac{1+q^2(y-x)-q^2(1-q)^2xy-\sqrt{(1+q^2(y-x)-q^2(1-q)^2xy)^2-4q^2y(1+q(1-q)x)^2}}
{2q^2y(1+q(1-q)x)},
\end{align}
which is the root of the kernel $1-\frac{qx(1-q+qz)(1+q(1-q)yz)}{(1-z)(1-q^2yz)}=0$. Then, by substitution of \eqref{z0ker1} into \eqref{eqCu1}, we get
\begin{align*}
J^u(x,y,q;1)=\frac{q^3y(z_0-1)}{(1+q(1-q)yz_0)}
-\frac{q^2yz_0(z_0-1)}{1-q^2yz_0}J^{bu}(x,y,q;z_0).
\end{align*}
We finally insert this expression into \eqref{eqCu1}, apply Lemma \ref{lem1}, and solve for $J^u(x,y,q;z)$. Hence,
\begin{lemma}\label{lem2_outer}
	The generating function $J^{u}(x,y,q;z)$ is given by
	\begin{align*}
	&\frac{q^4xy(2(1-z)\alpha^2+q(1-q^2yz)\beta
		(1-2q(1-z)-q\beta)x-q^3(1-q)\beta^3x^2)}
	{2((1-z)\alpha+q(1-q+qz)\beta x)(\alpha^2-q^2\beta^2x)}\\
	&+\frac{q^5x^2y(1+q+q^2(1-q)y-q^2(q-1)(1-(1-q)^2y)x)\beta}{
		2((1-z)\alpha+q(1-q+qz)\beta x)\sqrt{(1+q^2(y-x)-q^2(1-q)^2xy)^2-4q^2y\beta^2}},
	\end{align*}
	where $\alpha=1-q^2yz$ and $\beta=1+q(1-q)yz$.
\end{lemma}

\subsection{Enumeration over $\CP$}  We now lift up all the collected results to $\CP.$ To that goal, let $\nu\in\CP_k$. The polyomino $\nu$ satisfies exclusively of the following cases (see Figure \ref{figd3}):
\begin{figure}[htp]
	\begin{picture}(100,73)
	\setlength{\unitlength}{0.1cm}\put(0,30){
		%---1
		\put(0,15){
			\put(0,26){\tiny  Case (1)} \put(0,0){\line(1,0){4}}\put(4,0){\line(0,1){20}}\put(4,20){\line(-1,0){4}}\put(0,20){\line(0,-1){20}}
			\put(-2,-3){\tiny$y=-1/2$}\put(-2,21){\tiny$y=k-1/2$}}
		%---2
		\put(20,15){
			\put(0,26){\tiny  Case (2)} \put(0,0){\line(1,0){4}}\put(4,0){\line(0,1){20}}\put(4,20){\line(-1,0){4}}\put(0,20){\line(0,-1){20}}
			\put(4,7){\put(0,0){\line(1,0){4}}\put(4,0){\line(0,1){8}}\put(4,8){\line(-1,0){4}}\put(0,8){\line(0,-1){8}}}           \put(-2,-3){\tiny$y=-1/2$}\put(-2,21){\tiny$y=k-1/2$}
			\multiput(4,16)(0.8,0){12}{\circle*{.01}}
			\multiput(4,4)(1.5,0){7}{\line(1,0){.7}}
			\put(6,2){\tiny$y=1/2$}\put(-2,-3){\tiny$y=-1/2$}\put(-2,21){\tiny$y=k-1/2$}
			\put(4,17){\tiny$y=k-3/2$}}
		%---3
		\put(45,15){
			\put(0,26){\tiny  Case (3)} \put(0,0){\line(1,0){4}}\put(4,0){\line(0,1){20}}\put(4,20){\line(-1,0){4}}\put(0,20){\line(0,-1){20}}
			\put(4,0){\put(0,0){\line(1,0){4}}\put(4,0){\line(0,1){11}}\put(4,11){\line(-1,0){4}}\put(0,11){\line(0,-1){11}}}           \put(-2,-3){\tiny$y=-1/2$}\put(-2,21){\tiny$y=k-1/2$}
			\multiput(4,16)(0.8,0){12}{\circle*{.01}}
			\put(-2,-3){\tiny$y=-1/2$}\put(-2,21){\tiny$y=k-1/2$}
			\put(4,17){\tiny$y=k-3/2$}}
		%---4
		\put(70,15){
			\put(0,26){\tiny  Case (4)} \put(0,0){\line(1,0){4}}\put(4,0){\line(0,1){20}}\put(4,20){\line(-1,0){4}}\put(0,20){\line(0,-1){20}}
			\put(4,9){\put(0,0){\line(1,0){4}}\put(4,0){\line(0,1){11}}\put(4,11){\line(-1,0){4}}\put(0,11){\line(0,-1){11}}}           \put(-2,-3){\tiny$y=-1/2$}\put(-2,21){\tiny$y=k-1/2$}
			\multiput(4,4)(1.5,0){7}{\line(1,0){.7}}
			\put(6,2){\tiny$y=1/2$}}
		%---5
		\put(95,15){
			\put(0,26){\tiny  Case (5)} \put(0,0){\line(1,0){4}}\put(4,0){\line(0,1){20}}\put(4,20){\line(-1,0){4}}\put(0,20){\line(0,-1){20}}
			\put(4,0){\put(0,0){\line(1,0){4}}\put(4,0){\line(0,1){20}}\put(4,20){\line(-1,0){4}}\put(0,20){\line(0,-1){20}}}           \put(-2,-3){\tiny$y=-1/2$}\put(-2,21){\tiny$y=k-1/2$}}
	}
	%----6
	\put(0,10){
		\put(0,26){\tiny  Case (6)} \put(0,0){\line(1,0){4}}\put(4,0){\line(0,1){20}}\put(4,20){\line(-1,0){4}}\put(0,20){\line(0,-1){20}}
		\put(4,-8){\put(0,0){\line(1,0){4}}\put(4,0){\line(0,1){19}}\put(4,19){\line(-1,0){4}}\put(0,19){\line(0,-1){19}}}
		\multiput(4,16)(0.8,0){12}{\circle*{.01}}
		\put(-2,21){\tiny$y=k-1/2$}
		\put(-3,-3){\tiny$y=-1/2$}\put(4,17){\tiny$y=k-3/2$}}
	%----61
	\put(28,0){
		\put(0,36){\tiny  Case (6')}    \put(0,0){\line(1,0){4}}\put(4,0){\line(0,1){20}}\put(4,20){\line(-1,0){4}}\put(0,20){\line(0,-1){20}}
		\put(4,12){\put(0,0){\line(1,0){4}}\put(4,0){\line(0,1){19}}\put(4,19){\line(-1,0){4}}\put(0,19){\line(0,-1){19}}}
		\multiput(4,4)(1.5,0){7}{\line(1,0){.7}}
		\put(-1,21){\tiny$y=k-1/2$}
		\put(-3,-3){\tiny$y=-1/2$}\put(4,2){\tiny$y=1/2$}}
	%----7
	\put(50,10){
		\put(0,26){\tiny  Case (7)} \put(0,0){\line(1,0){4}}\put(4,0){\line(0,1){20}}\put(4,20){\line(-1,0){4}}\put(0,20){\line(0,-1){20}}
		\put(4,-8){\put(0,0){\line(1,0){4}}\put(4,0){\line(0,1){28}}\put(4,28){\line(-1,0){4}}\put(0,28){\line(0,-1){28}}}
		\put(-3,-3){\tiny$y=-1/2$}\put(-2,21){\tiny$y=k-1/2$}}
	%----71
	\put(70,0){
		\put(0,36){\tiny  Case (7')}    \put(0,0){\line(1,0){4}}\put(4,0){\line(0,1){20}}\put(4,20){\line(-1,0){4}}\put(0,20){\line(0,-1){20}}
		\put(4,0){\put(0,0){\line(1,0){4}}\put(4,0){\line(0,1){28}}\put(4,28){\line(-1,0){4}}\put(0,28){\line(0,-1){28}}}
		\put(-3,-3){\tiny$y=-1/2$}\put(-1,21){\tiny$y=k-1/2$}}
	%----8
	\put(90,5){
		\put(0,31){\tiny  Case (8)} \put(0,0){\line(1,0){4}}\put(4,0){\line(0,1){20}}\put(4,20){\line(-1,0){4}}\put(0,20){\line(0,-1){20}}
		\put(4,-4){\put(0,0){\line(1,0){4}}\put(4,0){\line(0,1){28}}\put(4,28){\line(-1,0){4}}\put(0,28){\line(0,-1){28}}}
		\put(-3,-3){\tiny$y=-1/2$}\put(-1,21){\tiny$y=k-1/2$}}
	\end{picture}
	\caption{Decomposition of a polyomino in $\CP_k$}\label{figd3}
\end{figure}
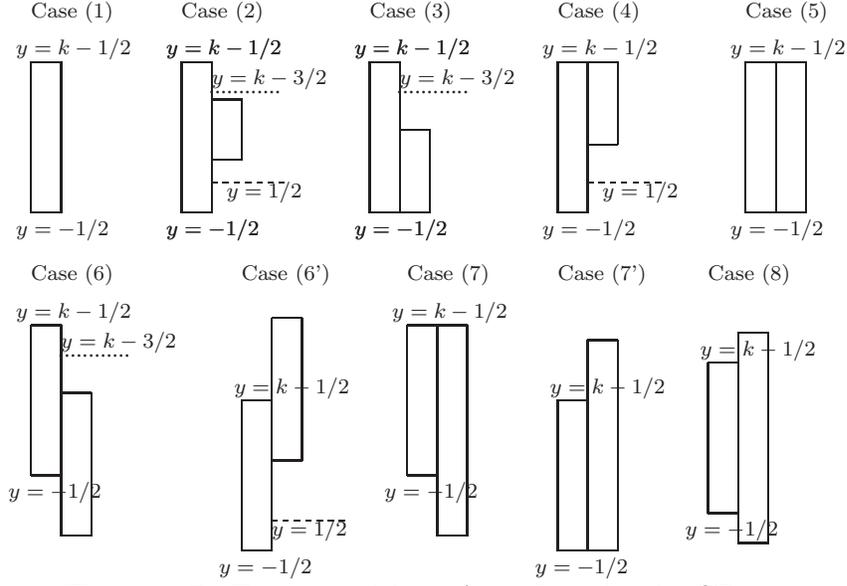
\begin{itemize}
	\item[Case (1)] It has one column;
	\item[Case (2)] It has at least two columns such that $u(\nu,2)<u(\nu,1)$ and $b(\nu,2)>b(\nu,1)$;
	\item[Case (3)] It has at least two columns such that $u(\nu,2)<u(\nu,1)$ and $b(\nu,2)=b(\nu,1)$;
	\item[Case (4)] It has at least two columns such that $u(\nu,2)=u(\nu,1)$ and $b(\nu,2)>b(\nu,1)$;
	\item[Case (5)] It has at least two columns such that $u(\nu,2)=u(\nu,1)$ and $b(\nu,2)=b(\nu,1)$;
	\item[Case (6)] It has at least two columns such that $u(\nu,2)<u(\nu,1)$ and $b(\nu,2)<b(\nu,1)$;
	\item[Case (6')] It has at least two columns such that $u(\nu,2)>u(\nu,1)$ and $b(\nu,2)>b(\nu,1)$;
	\item[Case (7)] It has at least two columns such that $u(\nu,2)=u(\nu,1)$ and $b(\nu,2)<b(\nu,1)$;
	\item[Case (7')] It has at least two columns such that $u(\nu,2)>u(\nu,1)$ and $b(\nu,2)=b(\nu,1)$;
	\item[Case (8)] It has at least two columns such that $u(\nu,2)>u(\nu,1)$ and $b(\nu,2)<b(\nu,1)$.
\end{itemize}
Next we write a recurrence for $J_k(x,y,q)$. This is done by adding up the contribution of each case to $J_k(x,y,q)$; that is,
\beq
J_k(x,y,q)&=&\underbrace{q^{2k+2}xy^k}_{(1)}
+\underbrace{\sum_{l=1}^{k-1}(k-1-l)xy^{k-l}q^{2k-2l}J^{bu}_l(x,y,q)}_{(2)}\\
&&+\underbrace{2\sum_{l=1}^{k-1}xy^{k-l}q^{2k-2l+1}J^u_l(x,y,q)}_{(3)+(4)}
+\underbrace{q^2x J_k(x,y,q)}_{(5)}\\
&&+\underbrace{2\sum_{l=2}^k\frac{xy^{k-l+1}q^{2k-2l+2}-xy^kq^{2k}}{1-q^2y}J_l^u(x,y,q)
	+2\sum_{l\geq k+1}\frac{xyq^2-xy^kq^{2k}}{1-q^2y}J_l^u(x,y,q)}_{(6)+(6')}\\
&&+2\underbrace{\sum_{l\geq k+1}xqJ_l(x,y,q)}_{(7)+(7')}+
\underbrace{\sum_{l\geq k+2}xq^2(l-1-k)J_l(x,y,q)}_{(8)}.
\feq
Multiplying by $z^{k-1}$, summing over $k\geq1$, and exchanging the order of double sums in each term, we obtain
\beq
J(x,y,q;z)&=&\frac{q^4xy}{1-q^2yz}+\frac{q^4xy^2z^2}{(1-q^2yz)^2}J^{bu}(x,y,q;z)+\frac{2q^3xyz}{1-q^2yz}J^u(x,y,q;z)\\
&&+q^2xJ(x,y,q;z)-\frac{2q^2xyz}{(1-z)(1-q^2yz)}(J^u(x,y,q;z)-J^u(x,y,q;1))\\
&&+\frac{qx(q-2+2z)}{(1-z)^2}(J(x,y,q;z)-J(x,y,q;1))
+\frac{q^2x}{1-z}\frac{\partial}{\partial z}J(x,y,q;z)\mid_{z=1}.
\feq
Grouping all the terms with $J(x,y,q;z)$ in one side, we then arrive at
\beq
&&\left(1-q^2x-\frac{qx(q-2+2z)}{(1-z)^2}\right)J(x,y,q;z)\\
&&\qquad=\frac{q^4xy}{1-q^2yz}+\frac{q^4xy^2z^2}{(1-q^2yz)^2}J^{bu}(x,y,q;z)
+\frac{2q^3xyz}{1-q^2yz}J^u(x,y,q;z)\\
&&\qquad-\frac{2q^2xyz}{(1-z)(1-q^2yz)}(J^u(x,y,q;z)-J^u(x,y,q;1))\\
&&\qquad-\frac{qx(q-2+2z)}{(1-z)^2}J(x,y,q;1)
+\frac{q^2x}{1-z}\frac{\partial}{\partial z}J(x,y,q;z)\mid_{z=1}.
\feq
Once again we leverage Kernel method to solve this equation. To that end, let $z$ take the values
\begin{align} \label{zker2}
z_+&=\frac{1+q(1-q)x+q\sqrt{(1-q^2)x^2+x}}{1-q^2x}\notag ,\\
z_-&=\frac{1+q(1-q)x-q\sqrt{(1-q^2)x^2+x}}{1-q^2x},
\end{align}
which are the roots of $1-q^2x-\frac{qx(q-2+2z)}{(1-z)^2}=0$. Then, by substitution of \eqref{zker2} into \eqref{eqCu1}, we obtain a system of two equations with two variables $J(x,y,q;1)$ and
$\frac{\partial}{\partial z}J(x,y,q;z)\mid_{z=1}$. We then solve this system to get
\begin{theorem}
	The generating function $J(x,y,q;1)$ for the number of convex polyominoes $\nu$ according to $v(\nu)$, $h(\nu)$ and $o(\nu)$ is given by
	\begin{align*}
	&J(x,y,q;1)=\frac{1}{2}q^2y(1-q^2y)\\
	&+\frac{z_+(1-q+qz_+)(1-q^2yz_-)}{z_+-z_-}J^{u}(x,y,q;z_+)
	-\frac{z_-(1-q+qz_-)(1-q^2yz_+)}{z_+-z_-}J^{u}(x,y,q;z_-)\\
	&-\frac{q^2yz_+^2(1-q^2yz_-)(1-z_+)}{2(z_+-z_-)(1-q^2yz_+)}J^{bu}(x,y,q;z_+)
	+\frac{q^2yz_-^2(1-q^2yz_+)(1-z_-)}{2(z_+-z_-)(1-q^2yz_-)}J^{bu}(x,y,q;z_-).
	\end{align*}
	where $J^{bu}(x,y,q;z)$ and $J^{u}(x,y,q;z)$ are given by Lemma \ref{lem1_outer} and Lemma \ref{lem2_outer}, respectively.
\end{theorem}

Note that $J(x,x,1;1)$ is exactly the generating function \eqref{cp_perim_gf} whose coefficient of $x^{2n}$ is \eqref{eqfcpf}. We next find the derivative of $J(x,x,q;1)$ at $q=1$ and obtain
\beq
\frac{\partial}{\partial q}J(x,x,q;1)\mid_{q=1}&=&
-\frac{2(36x^6-13x^5-156x^4+201x^3-98x^2+22x-2)x^2}{(1-4x)^3(1-2x)}\notag \\
&&+\frac{4x^4(x+2)(8x-3)}{\sqrt{1-4x}^5}.
\feq
This explicit form proves the results \eqref{tot_out_resp_perimeter} and \eqref{avg_out_resp_perimeter}.
Moreover, the generating function for the number of polyominoes
$\nu$ according to the perimeter and the outer-site perimeter is given by
$J(x^2,x^2,q;1)$. In particular, we state the following
result.
\begin{theorem}
	The generating function for the number of polyominoes $\nu$
	according to $o(\nu)$ is given by
	\begin{align*}
	&\frac{A(q)}{B(q)}-\frac{q^5(1+2q^2-q^3)(1+3q-2q^2)}
	{(1+q-q^2)(2-2q^3+q^4)(1-2q-2q^2)\sqrt{(1+q+q^2)(1-3q+q^2)}}\\
	&=q^4+2q^6+4q^7+10q^8+28q^9+77q^{10}+208q^{11}+586q^{12}+1572q^{13}+\cdots,
	\end{align*}
	where
	\begin{align*}
	A(q)&=q^4\bigl(2-q-18q^2-q^3+83q^4+51q^5-229q^6-250q^7+362q^8+597q^9-297q^{10}\\
	&\qquad-868q^{11}+124q^{12}+828q^{13}-48q^{14}-544q^{15}+55q^{16}+312q^{17}-200q^{18}\\
	&\qquad+48q^{19}-4q^{20}\bigr),\\
	B(q)&=(1+q-q^2)(1-2q-2q^2)(2-2q^3+q^4)(1+2q-4q^3-5q^4+3q^6+2q^7-2q^8)\\
	&\cdot(1-2q-4q^2+4q^3+11q^4-4q^5-13q^6+10q^7-2q^8).
	\end{align*}
\end{theorem}
By performing a singularity analysis (see \cite[Section VI]{FS} for a comprehensive review), we derive \eqref{outer_asympt}. Finally, we take the derivative of $J(x^2,x^2,q;1)$ at $x=1$ and obtain \eqref{avg_perim_res_outer}.

%---------------------------------------------------------------

\end{document}